\documentclass[11pt,a4paper]{article}

\usepackage{amsmath,amssymb}
\newcommand{\R}{\mathbb{R}}
\newcommand{\N}{\mathbb{N}}
\newcommand{\Z}{\mathbb{Z}}
\def\cA{{\mathcal A}}
\def\cC{{\mathcal C}}

\def\cS{{\mathcal S}}
\newcommand{\ee}{\varepsilon}
\renewcommand{\div}{{\rm div}\,}

\newcommand{\Frac}{\displaystyle \frac}

\newcommand{\Sum}{\displaystyle \sum}

\def\d{\partial}
\def\ddl{\dot \Delta_l}
\def\ddq{\dot \Delta_q}
\def\tilde{\widetilde}
\def\hat{\widehat}

\newcommand{\n}{\nabla}
\newcommand{\fd}{\frac{d}{2}}
\newcommand{\fdp}{\frac{d}{p}}
\newcommand{\p}{\partial}
\newcommand{\ql}{q_{l}}
\newcommand{\ul}{u_{l}}
\newcommand{\hl}{h_{l}}
\newcommand{\Fl}{F_{l}}
\newcommand{\Gl}{G_{l}}
\newcommand{\Rl}{R_{l}}
\newcommand{\Rlp}{R_{l}'}
\newcommand{\g}{\int_{\mathbb{R}^{d}}}

\newcommand{\qe}{q_\ee}
\newcommand{\ue}{u_\ee}
\newcommand{\dq}{\delta q}
\newcommand{\du}{\delta u}

\newtheorem{thm}{Theorem}
\newtheorem{lem}{Lemma}
\newtheorem{cor}{Corollary}
\newtheorem{prop}{Proposition}
\newtheorem{defi}{Definition}

\newtheorem{rem}{Remark}

\title{Convergence of capillary fluid models: from the non-local to the local Korteweg model}

\author{Fr\'ed\'eric Charve\footnote{Universit\'e Paris-Est Cr\'eteil, Laboratoire d'Analyse et de Math\'ematiques Appliqu\'ees (UMR 8050), 61 Avenue du G\'en\'eral de Gaulle, 94 010 Cr\'eteil Cedex (France). E-mail: frederic.charve@univ-paris12.fr}, Boris Haspot \thanks{Basque Center of Applied Mathematics, Bizkaia Technology Park, Building 500, 
E-48160, Derio (Spain), haspot@bcamath.org } \footnote{Karls Ruprecht Universit\"at Heidelberg, Institut for Applied Mathematics, Im Neuenheimer Feld 294,
D-69120 Heildelberg, Germany. Tel. 49(0)6221-54-6112}}

\date{}
\begin{document}

\maketitle

\begin{abstract} In this paper we are interested in the barotropic compressible Navier-Stokes system endowed with a non-local capillarity tensor depending on a small parameter $\ee$ such that it heuristically tends to the local Korteweg system. After giving some physical motivations related to the theory of non-classical shocks (see \cite{Lefloch}) we prove global well-posedness (in the whole space $\R^d$ with $d\geq 2$) for the non-local model and we also prove the convergence, as $\ee$ goes to zero, to the solution of the local Korteweg system.
\end{abstract}

\section{Introduction}
\subsection{Presentation of the models}

This introduction is widely inspired by the works of F.Coquel, D.Diehl, C.Merkle, and C.Rohde and we refer to \cite{Rohdehdr}, \cite{5CR} for a deep presentation of the capillary models.

The mathematical description of liquid-vapour phase interfaces has a long history and has been recently renewed in the 80's after the works of Dunn and Serrin (see \cite{3DS}). The first investigations begin with the Young-Laplace theory which claims that the phases are separated by a hypersurface and that the jump in the pressure across the hypersurface is proportional to the curvature of the hypersurface. The main difficulty consists in describing the location and the movement of the interfaces.\\

Another major problem is to understand whether the interface behaves as a discontinuity in the state space (sharp interface, SI) or whether the phase boundary corresponds to a more regular transition (diffuse interface, DI).
The DI ansatz is often favored as it covers topological changes in the phase distribution (such as the separation of bubbles/drops). Moreover, instead of one system per phase coupled with free-boundary problems required by the SI model, in the DI approach only one set of equations has to be solved in a single spatial domain (the density takes into account the different phases) which considerably simplifies the mathematical and numerical study.\\

Another approach corresponds to determine equilibrium solutions which classically consists in the minimization of the free energy functional.
Unfortunately this minimization problem has an infinity of solutions, and many of them are physically wrong (some details are given later). In order to overcome this difficulty, Van der Waals in the XIX-th century was the first to add a term related to the surface energy to select the physically correct solutions, modulo the introduction of a diffuse interface. This theory is widely accepted as a thermodynamically consistent model for equilibria.\\

The study of the dynamics of a liquid-vapour mixture is complicated: the model is based upon the compressible Navier-Stokes equations with a Van der Waals state law for ideal fluids. Endowed with such a state law, the compressible Navier-Stokes equation can describe the dynamics of multi-phases (liquid and vapour) but the effect of surface tension is still not taken into account. In this goal, Korteweg extended (in \cite{3K}) the Navier-Stokes equations by adding a capillarity tensor modeling the behavior at the interfaces. The additional term is directly related to the surface energy added to the free functional in order to obtain physically relevant minimizers. This local Navier-Stokes-Korteweg system belongs to the class of diffuse interface models and, to penalize the density fluctuations, it involves the density gradient and in fact it introduces, in the classical model, third-order derivatives of the density (and therefore numerical complications).\\

Alternatively, another way to penalize the high density variations consists in applying a zero order but non-local operator to the density gradient \cite{9Ro}, \cite{5Ro}. Let us mention in passing that it was the original idea of Van der Waals in \cite{VW} and that statistical mechanics tend to show that this non-local approach is in fact the correct one (see \cite{Rohdehdr}). This model does not introduce additional orders of differentiation but a non-local integral operator which creates comparable numerical difficulties. For new developments  we also refer to the work of M. Heida and J. M\'alek in \cite{HM}.\\
In this introduction we will briefly recall both of these approaches and explain their main differences.

\subsection{Basic modelization setting}
We want here to explain and give the main lines of how were modeled the local and non-local Navier-Stokes-Korteweg systems. Let us consider a single-component fluid cointained in an open set $\Omega\in\R^d$ at constant temperature $T^{*}$. In the isothermal setting the density is denoted by $\rho: \Omega\rightarrow ]0,b[$ ($b>0$), and the pressure is given by $p:\rho\in]0,b[\mapsto p(\rho)\in]0,+\infty[$. To denote the two phases, we need to define the energy density function as:
\begin{equation}
W(\rho)=\rho f(\rho,T^{*}),\quad \rho\in]0,b[,\mbox{ and }f\in C^{2}(0,b). 
\label{2.1}
\end{equation}
From the mathematical point of view, $W$ has to satisfy the following properties:
\begin{itemize}
\item $\exists \alpha_{1},\alpha_{2}\in]0,b[$: $\alpha_{1}<\alpha_{2}$, $W^{''}>0$ in $]0,\alpha_{1}[\cup ]\alpha_{2},b[$, $W^{''}<0$ in $]\alpha_{1},\alpha_{2}[$.
\item $\lim_{\rho\rightarrow 0}W(\rho)= \lim_{\rho\rightarrow b}W(\rho)=+\infty$.
\item $W\geq 0$ in $]0,b[$.
\end{itemize}
We refer to \cite{9Ro} and \cite{5CR} for a far more complete presentation (and also figures). Since we consider the fluid as a single system, the density will play the role of an order parameter in the sense that it will determine the phase state. More precisely the high variation of density will define the location of the interfaces and we will say that the fluid at the position $x$ is in the \textit{vapour/elliptic (or spinodal)/liquid} phase if $\rho(x)\in ]0,\alpha_1]/[\alpha_1,\alpha_2]/[\alpha_2,b[$.\\
According standard thermodynamical theory, we define the pressure as follows:
\begin{equation}
P(\rho)=\rho W'(\rho)-W(\rho).
\label{pression}
\end{equation}
Thus the pressure is a {\it{non-convex}} and {\it{non-monotone}} function ( because $P^{'}(\rho)=\rho W^{''}(\rho)$). In passing, let us mention that the free energy $f$ cannot be directly measured while the pressure can be determined from experiments. A standard model for the pressure in two-phase mixtures is the well-known Van der Waals law given by:
\begin{equation}
P(\rho,T^{*})=\frac{RT^*\rho}{b-\rho}-a\rho^2
\label{VDW}
\end{equation}
where $R$ is the specific gas constant and $a$, $b$ are positive constants. Let us briefly say that the parameter $a$ controls the attractive forces between molecules of the fluid and $b$ is related to the molecule size.
\subsection{Variational Approach for sharp interfaces}
In the sequel, we want to model our fluid by considering solutions as static equilibria at constant temperature $T^*$ (at least in an asymptotic configuration in time). In this optic, we can associate to this problem of equilibrium for the mixture of liquid-vapor a variational problem that consists in minimizing the free energy functional. For this, we introduce an admissibility set for the density which takes into account the fact that our interfaces have null thickness, more precisely we set:
$$A^0=\{\rho\in L^{1}(\Omega)/   W(\rho)\in L^{1}(\Omega),\int_\Omega\rho(x)dx =m\}.$$
Here we have prescribed the total mass by a constant $m>0$. To find a static equilibrium, we need to minimize, when $\rho \in A^0$, the functional:
\begin{equation}
F^0[\rho ]=\int_\Omega W(\rho(x))dx.
\label{fonctionel}
\end{equation}
The problem is that it is possible to build an infinite number of minimizers: if $\beta_1\in]0,\alpha_1[$ and $\beta_2\in]\alpha_2,\beta[$ are the Maxwell states (the two points where $W$ is minimal), up to a translation we can have $W(\beta_1)=W(\beta_2)=0$ and then if $\beta_1|\Omega|<m<\beta_2|\Omega|$ every piecewise constant function $\rho$ taking its values in $\{\beta_1, \beta_2\}$, so that the total mass is $m$, is a minimizer. This allows drops and bubbles (which are physically correct solutions: they only take the values $\beta_1$ and $\beta_2$, and away from the boundary, the length of the phase interface is minimal) as well as an infinity of minimizers describing an arbitrary large number of phase changes (provided that the total mass is $m$) which clearly are physically wrong.\\

In other words, we need selection principles to identify the physically relevant solutions.

\begin{rem}
 \sl{All the previous minimizers realize the phase transition as a jump between the vapour or liquid zones, without taking values in the spinodal zone (at least up to set of zero Lebesgue measure). This is why we talk about sharp-interface model.}
\end{rem}


\subsection{The local diffuse interface approach}
Van der Waals seems to be the first who tried to overcome the lack of uniqueness for the sharp interface approach by adding a term of capillarity (see \cite{VW}). His idea for selecting some relevant physical solution is to penalize the high variations of density (typically at the interfaces) by adding some derivative terms to the previous functional. More precisely, he considers the admissibility set $A^{local}$ by:
$$A^{local}=H^1(\Omega)\cap A^0.$$
We now search a function $\rho_{\ee}$ which is a minimizer of the local Van der Waals functional:
$$F^{\ee}_{local}[\rho_\ee]=\int_{\Omega}(W(\rho_{\ee}(x))+\gamma\frac{\ee^2}{2}|\n\rho_{\ee}(x)|^2)dx,$$
where $\gamma>0$ is a capillarity coefficient and $\ee>0$ is a scaling parameter. Obviously, functions which exhibit jumps across a surface are not in $H^1(\Omega)$. By consequence, a possible two-phase minimizer $\rho_{\ee}$ cannot change of phase without passing by values in the elliptic region, another way to express this phenomena is to say that the constituent propagates continuously in the interfaces. This is the reason why the variational problem is called a diffuse interface approach. Here the thickness of the interface is expected to be of size $O(\ee)$.
\\

This diffuse interface approach may appear artificial at first sight but a result, obtained by Modica in \cite{Mo}, validates this method. He proves that if a sequence of minimizers $\rho_{\ee}$ converges when $\ee$ tends to zero, the limit is a physically relevant minimizer of the sharp-interface functional $F^0$.

\subsection{The non-local diffuse interface approach}
This alternative approach (also called global diffuse interface approach, introduced by Serrin et al in \cite{3DS} and next by Coquel, et al. in \cite{5CR}, and Rohde in \cite{5Ro} and \cite{Rohdehdr}) also consists in a modification of the sharp-interface functional $F_{0}$. One of the main interest of this approach is that it does not need to introduce high derivatives on the density to penalize the high fluctuations of density.\\
Let us consider the following admissibility set:
$$A^{global}=A^{0}\cap L^2(\Omega),$$
Let us choose a function $\phi\in L^1(\R^d)$ such that:
$$(|.|+|.|^2)\phi(.)\in L^1(\R^d)\mbox{, }\quad\int_{\R^d}\phi(x)dx=1,\quad\phi\mbox{ even, and }\phi\geq0.$$
Such a function is called an interaction potential, and $\phi_{\ee}(x)=\frac{1}{\ee^{d}}\phi(\frac{x}{\ee})$ is called scaled interaction potential. Then for $\gamma>0$, we search for $\rho_\ee\in A^{global}$ that minimizes the following modified functionnal called non-local Van der Waals functional:
$$F^{\ee}_{global}[\rho_{\ee}]=\int_{\Omega}\Big(W(\rho_{\ee}(x))+\frac{\gamma}{4}\int_{\Omega}
\phi_{\ee}(x-y)(\rho_{\ee}(y)-\rho_{\ee}(x))^2dy\Big)dx.$$
Roughly speaking, the non-local term penalizes high fluctuations on an $\ee$-scale if the interaction potential has most of its mass in a ball centered around zero and whose radius is of order $\ee$. The set $A^{global}$ contains functions with jumps but it can be proved that minimizers of $F_{global}^\ee$ take values also in the elliptic regions: the model belongs to the class of diffuse-interface models. The counterpart of the result of Modica was proved by Alberti and Bellettini (\cite{Alb}), which also validates this approach.

\begin{rem}
\sl{Let us consider $F^{\ee}_{global}$ with $\Omega=\R$ and assume that $\rho\in C^{\infty}(\R)$ is an analytic function with the expansion (we refer to \cite{Rohdehdr}):
$$\rho(y)=\rho(x)+\sum_{k\geq 1}\frac{1}{k!}(y-x)^{k}\frac{\p^{k}}{\p x^{k}}\rho(x),$$
If we put this expression into the non-local term in $F^{\ee}_{global}$ and perform the change of variable $r=\frac{x-y}{\ee}$ we obtain:
$$\int_{\R}\int_{\R}\phi(r)(\sum_{k\geq 1}\frac{1}{k!}(\ee r)^{k}\frac{\p^{k}}{\p x^{k}}\rho(x))^2 drdx\approx \int_{\R}\int_{\R}\ee^{2}\phi(r)(r\frac{\p}{\p x}\rho(x))^2 drdx.$$
Since $\int_{\R}\phi(r)r^{2}dr<+\infty$ holds by definition of $\phi$, the last expression is equal to the local penalty terms in $F^{\ee}_{local}$ (up to a multiplicative constant). In particular the scaling with respect to the small parameter $\ee$ fits.}
\end{rem}

\subsection{Non classical shocks and the Euler problem}
Another important approach using the capillarity concerns the discontinuous solutions of nonlinear hyperbolic systems or conservation laws. Indeed the mathematical model of liquid-vapor flows should have special solutions that can be interpreted as dynamical phase transitions. In \cite{Rohdehdr}, the author shows that there are traveling-wave solutions for one-dimensional versions of the local and non local NSK system that connects states in different phase. This fact is one of the major arguments to accept the $NSK$ systems as promising candidates to model the dynamics of liquid vapour flow in a reliable way. 

We also refer to the works of \cite{BDDJ}  in the case of the local diffuse approach where the authors prove the existence of planar traveling waves representing either diffuse interfaces or solitons. The system of equations governing planar traveling waves reduces to a planar Hamiltonian system, for which a portrait analysis exhibits heteroclinic/homoclinic orbits (see for more details \cite{BDDJ}).\\

 The first order model that governs the sharp-interface limit is of particular interest since it is not a standard conservation law (we will explain in the sequel more precisely what it means). In particular it is interesting in this sense to consider the solutions of Euler systems as limit solutions of the NSK systems.
\\
In the sequel, we will interest in the Euler or sharp-interface case, we look for a specific volume $\tau$ in $(\frac{1}{b},\infty)$ and a velocity $v$ that satisfy in $(0,+\infty)\times\R$ the equations:
\begin{equation}
\begin{aligned}
&\p_{t}\tau-\p_{x}v=0,\\
&\p_{t}v-\p_{x}(\widetilde{P}(\tau))=0,
\end{aligned}
\label{euler}
\end{equation}
with the function $\widetilde{P}:(\frac{1}{b},\infty)\rightarrow (0,\infty)$ given by:
$$\widetilde{P}(\tau)=P(\frac{1}{\tau}),\;\;\;\tau\in(\frac{1}{b},\infty).$$
More precisely if we investigate the well-known Euler system with a Van der Waals pressure state in one dimension, we can recall that following the theory of conservation law; the two eigenvalues of the system are:
\begin{equation}
\lambda_{1}(\tau,v)=-\sqrt{-\widetilde{P}^{'}(\tau)},\;\;\;\lambda_{2}(\tau,v)=-\sqrt{-\widetilde{P}^{'}(\tau)}.
\label{vp}
\end{equation}
The corresponding eigenvectors $r_{1}$, $r_{2}$ are:
\begin{equation}
w_{1}(\tau,v)=\left(\begin{array}{c}
1\\
\sqrt{-\widetilde{P}^{'}(\tau)}\\
\end{array}
\right),\;\;w_{2}(\tau,v)=\left(\begin{array}{c}
1\\
-\sqrt{-\widetilde{P}^{'}(\tau)}\\
\end{array}
\right)
\end{equation}
Furthermore by calculus we obtain:
\begin{equation}
\n\lambda_{1}(\tau,v)\cdot w_{1}(\tau,v)=\frac{\widetilde{P}^{''}(\tau)}{2\sqrt{-\widetilde{P}^{'}(\tau)}},\;\;
\n\lambda_{2}(\tau,v)\cdot w_{2}(\tau,v)=\frac{-\widetilde{P}^{''}(\tau)}{2\sqrt{-\widetilde{P}^{'}(\tau)}}
\end{equation}
We now recall the definition of a \textit{standard conservation law} in the sense of Lax (it means entropy solutions):
\begin{itemize}
\item The system is \textbf{strictly entropic} if the eigenvalues are distinct and real.
\item The characteristics fields are \textbf{genuinely nonlinear} if we have for all $(\tau,v)$,
$$\n\lambda_{1}(\tau,v)\cdot w_{1}(\tau,v)\ne 0\;\;\mbox{and}\;\;\n\lambda_{2}(\tau,v)\cdot w_{2}(\tau,v)\ne 0,$$
\end{itemize}
for more details we refer to \cite{Serre}. The definition of genuine nonlinearity is some kind of extension of the notion of convexity to vector-valued functions (in particular when we consider the specific case of the traveling waves).  The previous assumptions aim at assuring the existence and the uniqueness of the Riemann problem ( see \cite{Evans} and \cite{Serre}).  When $P$ is a Van der Waals pressure, we observe that then the first conservation law (\cite{Serre,Evans}) is far from being a standard hyperbolic system:
\begin{itemize}
\item  It is not hyperbolic (but elliptic) in $(\frac{1}{\alpha_{1}},\frac{1}{\alpha_{2}})\times\R$,
\item the characteristic fields are not genuinely nonlinear in the hyperbolic part of the state space.
\end{itemize}
Here the classical theory cannot be applied. It has been shown (see \cite{Evans}) that there exists an infinity of weak \textbf{entropy} solutions of the Riemann problem with initial states in different phases. It means that the entropy inequalities are not sufficiently discriminating when one characteristic field fails to be genuinely nonlinear. \textbf{One of the main research in this field is to understand what are the physical relevant solutions when the flux is non genuinely non linear}. An idea for selecting some relevant solutions of the Euler system with a Van der Waals pressure is to consider the solutions limit of solutions of $NSK$ systems when the capillarity and the viscosity coefficients tends to $0$. This condition is called the viscosity-capillarity criterium. By a study of traveling waves (see \cite{Rohdehdr} or \cite{Lefloch} in the context of the Korteweg de Vries system) a remarkable fact is to observe that generally the limit solutions violate the Oleinik conditions, more precisely the shock are undercompressive (it means non classic in the sense of the Lax theory). We refer for more details to the book \cite{Lefloch}  of P. Lefloch.\\
An important research line (see \cite{5Ro}, \cite{Rohdehdr}, \cite{5CR}) is to model the capillarity tensor and to understand how fast the solutions converges to the Euler system
when the capillarity and the viscosity coefficients tends to zero. We want point out here that it exists three different regimes, more precisely if we assume the viscosity coefficient equal to $\ee$ with $\ee\rightarrow 0$. Then we have the three different regimes:
\begin{itemize}
\item $\kappa<<\ee^{2}$, the viscosity dominates.
\item $\kappa \simeq\ee^{2}$, intermediary regime.
\item $\kappa>>\ee^{2}$, the capillarity dominates.
\end{itemize}
A second important problem from a numerical point of view is to describe the link between the local and the non local Korteweg system and more particularly how it plays on the choice of the thickness of the interfaces. We propose us to answer to this second question by giving a speed of convergence between the local and non local Korteweg system after suitable choices of kernel for the local system.
\subsection{Results}
In this work we wish to unify the two models, the local Korteweg system and the non local Korteweg system, more precisely we will search to approach the local Korteweg solutions by solutions of the non local Korteweg system in an appropriate way (it means by choosing suitable
kernel on the non local tensor of capillarity).\\

Let $\rho$ and $u$ denote the density and the velocity of a compressible viscous fluid. As usual, $\rho$ is a non-negative function and $u$ is a vector-valued function defined on $\R^d$. In the sequel we will denote by $\cA$ the following diffusion operator
$$
\cA u= \mu \Delta u+ (\lambda+\mu)\nabla \div u, \quad \mbox{with} \quad \mu>0 \quad \mbox{and} \quad \nu=\lambda+ 2 \mu >0.
$$
The Navier-Stokes equation for compressible fluids endowed with internal non-local capillarity introduced in \cite{5Ro} reads:
$$
\begin{cases}
\begin{aligned}
&\d_t\rho+\div (\rho u)=0,\\
&\d_t (\rho u)+\div (\rho u\otimes u)-\cA u+\nabla(P(\rho))=\kappa\rho\nabla(\phi*\rho-\rho),\\
\end{aligned}
\end{cases}
\leqno{(NSR)}
$$
supplemented by the following conditions on the kernel $\phi$:
$$(|.|+|.|^2)\phi(.)\in L^1(\R^d)\mbox{, }\quad\int_{\R^d}\phi(x)dx=1,\quad\phi\mbox{ even, and }\phi\geq0.$$
The mass, momentum equations for the Korteweg system read ( see \cite{3DS}):
$$
\begin{cases}
\begin{aligned}
&\d_t\rho+\div (\rho u)=0,\\
&\d_t (\rho u)+\div (\rho u\otimes u)-\cA u+\nabla(P(\rho))=\kappa\rho\nabla\Delta\rho,\\
\end{aligned}
\end{cases}
\leqno{(NSK)}
$$
Both of these systems have been studied in the context of existence of strong solutions in critical spaces for the scaling of the equations. For example, concerning the strong solutions, we refer to \cite{DD} for a study of $(NSK)$ system, and to \cite{Has1} for $(NSR)$.

If we compute the Fourier transform in the capillarity terms (simplified by the density), we obtain $(\hat{\phi}(\xi)-1) \hat{\rho}(\xi)$ in the non-local model, and $-|\xi|^2 \hat{\rho}(\xi)$ in the local model. We want to see if, when formally $\hat{\phi}(\xi)$ is "close" to $1-|\xi|^2$, we can expect the solutions of these models to be close.

The aim of this paper is to approximate the local Korteweg model $(NSK)$ with a non-local model such as system $(NSR)$ introduced on his modern form by C. Rohde in \cite{5Ro} and also \cite{5CR} (see Van der Waaals \cite{VW} for the originally works) . For that we will choose a specific function $\phi_{\ee}$ in the capillarity tensor, and the system we will consider in this paper is the following:
$$
\begin{cases}
\begin{aligned}
&\d_t\rho_\ee+\div (\rho_\ee \ue)=0,\\
&\d_t (\rho_\ee \ue)+\div (\rho_\ee u\otimes \ue)-\cA \ue+\nabla(P(\rho_\ee))=\rho_\ee\frac{\kappa}{\ee^2}\nabla(\phi_\ee*\rho_\ee-\rho_\ee),\\
\end{aligned}
\end{cases}
\leqno{(NSR_\ee)}
$$
where we set:
$$\phi_{\ee}=\frac{1}{\ee^d} \phi(\frac{x}{\ee}) \quad \mbox{with} \quad \phi(x)=\frac{1}{(2 \pi)^d} e^{-\frac{|x|^2}{4}}$$
As the Fourier transform of $\phi$ is $\hat{\phi}(\xi)=e^{-|\xi|^2}$ we have: $$\hat{\phi_\ee}(\xi)=e^{-\ee^2|\xi|^2},$$
and for a fixed $\xi$, when $\ee$ is small, $\Frac{\hat{\phi_\ee}(\xi)-1}{\ee^2}$ is close to $-|\xi|^2$.
\begin{rem}
We mention here that the choice of $\phi_\ee$ is in accordance with the physical relevant capillarity coefficient (see \cite{Rohdehdr}).
\end{rem}
We will consider a density which is close to an equilibrium state $\overline{\rho}$ and we will introduce the change of function $\rho= \overline{\rho}(1+q)$.  By simplicity we take $\overline{\rho}=1$. The previous systems become:
$$
\begin{cases}
\begin{aligned}
&\d_t q+ u.\nabla q+ (1+q)\div u=0,\\
&\d_t u+ u.\nabla u -\cA u+P'(1).\nabla q-\kappa\nabla \Delta q= K(q).\nabla q- I(q) \cA u,\\
\end{aligned}
\end{cases}
\leqno{(K)}
$$
and
$$
\begin{cases}
\begin{aligned}
&\d_t \qe+ \ue.\nabla \qe+ (1+\qe)\div \ue=0,\\
&\d_t \ue+ \ue.\nabla \ue -\cA \ue+P'(1).\nabla \qe-\frac{\kappa}{\ee^2}\nabla(\phi_\ee*\qe-\qe)\\
&\hspace{7cm}= K(\qe).\nabla \qe- I(\qe) \cA \ue,\\
\end{aligned}
\end{cases}
\leqno{(R_\ee)}
$$
where $K$ and $I$ are real-valued functions defined on $\R$ given by:
$$
K(q)=\left(P'(1)-\frac{P'(1+q)}{1+q}\right) \quad \mbox{and} \quad I(q)=\frac{q}{q+1}.
$$
\begin{rem}
\sl{In the sequel we will sometimes rewrite the term $K(q).\nabla q$ as $\nabla (G(q))$ where $G$ is a primitive of $K$.}
\end{rem}

\begin{rem}
\sl{When $\overline{\rho}\neq1$ the only changes are in the viscosity and capillarity coefficients, and in the expression of functions $I$ and $K$.}
\end{rem}
We now want recall differents results concerning the local and non-local Korteweg system by emphasizing on the results of global weak solutions.
\subsubsection*{Results on the global weak solutions for the local Korteweg system}
The existence of global weak solution
for the model of Korteweg with constant capillary coefficient for $N\geq 2$
is still an open problem. Indeed the main difficulty states in treating the capillarity tensor and more specially the quadratic terms in gradient of the density, as a matter of fact it seems impossible to pass to the limit in this term as $|\n\rho|^{2}$ is only in $L^\infty(L^1)$. We can only hope a convergence in the sense of the measure. A way to overcome this difficulty would be to obtain regularizing effects on the density. In this spirit D. Bresch, B. Desjardins and C-K. Lin in \cite{3BDL} got some global weak solutions for the isotherm Korteweg model with some
specific viscosity coefficients $\mu(\rho)=C\rho$ with $C>0$ and $\lambda(\rho)=0$. By choosing
these specific coefficients they obtain a gain of derivatives on the
density $\rho$ where $\rho$ belongs to $L^{2}(H^{2})$. Recently the second author in \cite{Has3}has shown the existence of global weak solution with small initial data in the energy space for specific choices of the capillary coefficients ($\kappa(\rho)=\frac{K}{\rho^2}$ with $K>0$) and with general viscosity coefficient.
Comparing with the results of \cite{3BDL}, he gets global weak solutions with test function $\varphi\in C^{0}_{\infty}(\R^d)$  independant of the density $\rho$ as it is the case in \cite{3BDL}.
In fact he is able to obtain a gain of derivative on the density $\rho$ by using the specific structure of the capillary tensor. In particular it is related with the Gross-Pitaevski system via the Madelung transformation (see \cite{Has3} for more details).
\subsubsection*{Results on the existence of strong solution for the local Korteweg system}
Let us briefly mention that the existence of strong solutions for $N\geq2$ is known since the works by H. Hattori and D. Li \cite{3H1,3H2}. R. Danchin and B. Desjardins in \cite{DD} improve this result by working in critical spaces for the scaling of the equations, more precisely the initial data $(\rho_{0},\rho_{0}u_{0})$
belong to $B^{\frac{d}{2}}_{2,1}\times B^{\frac{d}{2}1}_{2,1}$ (the fact that $B^{\frac{d}{2}}_{2,1}$ is embedded in $L^{\infty}$ plays a crucial role to control the vacuum but also for some reason concerning the multiplier spaces). In \cite{3MK}, M. Kotschote showed the existence of strong solution for the isothermal model in bounded domain by using Dore\^a-Venni Theory and $\mathcal{H}^{\infty}$ calculus.
In \cite{Has2},  the second author generalizes the results of \cite{DD} in the case of non isothermal Korteweg system with physical coefficients depending on the density and the temperature. He gets strong solutions with initial data belonging to the critical spaces $B^{\frac{d}{2}}_{2,1}\times B^{\frac{d}{2}-1}_{2,1}\times B^{\frac{d}{2}-2}_{2,1}$ when the physical coefficients depend only on the density.
\subsubsection*{Existence of global weak solution and strong solutions for the non-local Korteweg system}
The existence of global weak solutions was obtained by the second author in \cite{Has4} by following the Lion's theory. In passing, we would mention that the capillarity allows to treat general pressure of the form $P(\rho)=a\rho^\gamma$ with $\gamma>1$, indeed we get automatically some estimate on $\rho$ in $L^2_{loc}$ what allows us to get renormalized solutions without efforts.
The first result of well-posedness in finite time are coming from the works of C. Rohde in \cite{5Ro}. More recently the second authors in \cite{Has1} improved these results by obtaining strong solutions in initial data invariant for the scaling of the system.

\subsection{Littlewood-Paley theory and Besov spaces}

\subsubsection{Littlewood-Paley theory}

In the sequel the Fourier transform of $u$ with respect to the space variable will be denoted by $\mathcal{F}(u)$ or $\hat{u}$. 
In this section we will state classical definitions and properties concerning the homogeneous dyadic decomposition with respect to the Fourier variable. We will recall isome classical results and we refer to \cite{Dbook} (Chapter 2) for proofs (and more general properties), and to the appendix for specific properties used in this paper.

To build the Littlewood-Paley decomposition, we need to fix a smooth radial function $\chi$ supported in (for example) the ball $B(0,\frac{4}{3})$, equal to 1 in a neighborhood of $B(0,\frac{3}{4})$ and such that $r\mapsto \chi(r.e_r)$ is nondecreasing over $\R_+$. So that if we define $\varphi(\xi)=\chi(\xi/2)-\chi(\xi)$, then $\varphi$ is compactly supported in the annulus $\{\xi\in \R^d, \frac{3}{4}\leq |\xi|\leq \frac{8}{3}\}$ and we have that,
\begin{equation}
 \forall \xi\in \R^d\setminus\{0\}, \quad \sum_{j\in\Z} \varphi(2^{-j}\xi)=1.
\label{LPxi}
\end{equation}
Then we can define the \textit{dyadic blocks} $(\ddq)_{q\in \Z}$ by $\ddq:= \varphi(2^{-q}D)$ (that is $\hat{\ddq u}=\varphi(2^{-q}\xi)\hat{u}(\xi)$) so that, formally, we have
\begin{equation}
u=\Sum_q \ddq u
\label{LPsomme} 
\end{equation}
As (\ref{LPxi}) is satisfied for $\xi\neq 0$, the previous formal equality holds true for tempered distributions \textit{modulo polynomials}. A way to avoid working modulo polynomials is to consider the set $\cS_h'$ of tempered distributions $u$ such that
$$
\lim_{q\rightarrow -\infty} \|\dot{S}_q u\|_{L^\infty}=0,
$$
where $\dot{S}_q$ stands for the low frequency cut-off defined by $\dot{S}_q:= \chi(2^{-q}D)$. If $u\in \cS_h'$, (\ref{LPsomme}) is true and we can write that $\dot{S}_q u=\Sum_{p\leq q-1} \ddq u$. We can now define the homogeneous Besov spaces used in this article:
\begin{defi}
\label{LPbesov}
 For $s\in\R$ and  
$1\leq p,r\leq\infty,$ we set
$$
\|u\|_{\dot B^s_{p,r}}:=\bigg(\sum_{q} 2^{rqs}
\|\Delta_q  u\|^r_{L^p}\bigg)^{\frac{1}{r}}\ \text{ if }\ r<\infty
\quad\text{and}\quad
\|u\|_{\dot B^s_{p,\infty}}:=\sup_{q} 2^{qs}
\|\Delta_q  u\|_{L^p}.
$$
We then define the space $\dot B^s_{p,r}$ as the subset of  distributions $u\in {\cS}'_h$ such that $\|u\|_{\dot B^s_{p,r}}$ is finite.
\end{defi}
Once more, we refer to \cite{Dbook} (chapter $2$) for properties of the inhomogeneous and homogeneous Besov spaces. Among these properties, let us mention:
\begin{itemize}
\item for any $p\in[1,\infty]$ we have the following chain of continuous embeddings:
$$
\dot B^0_{p,1}\hookrightarrow L^p\hookrightarrow \dot B^0_{p,\infty};
$$
\item if $p<\infty$ then 
  $B^{\frac dp}_{p,1}$ is an algebra continuously embedded in the set of continuous 
  functions decaying to $0$ at infinity;
  \item the following real interpolation property is satisfied for
  all $1\leq p,r_1,r_2,r\leq\infty,$ $s_1\not=s_2$ and $\theta\in(0,1)$:
  $$
  [\dot B^{s_1}_{p,r_1},\dot B^{s_2}_{p,r_2}]_{(\theta,r)}=\dot B^{\theta s_2+(1-\theta)s_1}_{p,r};
  $$
  \item for any  smooth homogeneous  of degree $m$ function $F$ on $\R^d\setminus\{0\}$
the operator $F(D)$ maps  $\dot B^s_{p,r}$ in $\dot B^{s-m}_{p,r}.$ This implies that the gradient operator maps $\dot B^s_{p,r}$ in $\dot B^{s-1}_{p,r}.$  
  \end{itemize}
The following lemma (referred as \emph{Bernstein's inequalities}) describes the way how derivatives act on spectrally localized functions.
  \begin{lem}
\label{LPbernstein}
{\sl
Let  $0<r<R.$   A
constant~$C$ exists so that, for any nonnegative integer~$k$, any couple~$(p,q)$ 
in~$[1,\infty]^2$ with  $q\geq p\geq 1$ 
and any function $u$ of~$L^p$,  we  have for all $\lambda>0,$
$$
\displaylines{
{\rm Supp}\, \widehat u \subset   B(0,\lambda R)
\Longrightarrow
\|D^k u\|_{L^q} \leq
 C^{k+1}\lambda^{k+d(\frac{1}{p}-\frac{1}{q})}\|u\|_{L^p};\cr
{\rm Supp}\, \widehat u \subset \{\xi\in\R^N\,/\, r\lambda\leq|\xi|\leq R\lambda\}
\Longrightarrow C^{-k-1}\lambda^k\|u\|_{L^p}
\leq
\|D^k u\|_{L^p}
\leq
C^{k+1}  \lambda^k\|u\|_{L^p}.
}$$
}
\end{lem}   
 The first Bernstein inequality entails the following embedding result:
\begin{prop}\label{LP:embed}
\sl{For all $s\in\R,$ $1\leq p_1\leq p_2\leq\infty$ and $1\leq r_1\leq r_2\leq\infty,$
  the space $\dot B^{s}_{p_1,r_1}$ is continuously embedded in 
  the space $\dot B^{s-d(\frac1{p_1}-\frac1{p_2})}_{p_2,r_2}.$}
\end{prop}

In this paper, we shall mainly work with functions or distributions depending on both the time variable $t$ and the space variable $x.$ Most often, these functions will be seen as defined on some time interval $I$ and  valued in some Banach space $X.$ We shall denote by $\cC(I;X)$ the set of continuous functions on $I$ with values in $X.$ For $p\in[1,\infty]$, the notation $L^p(I;X)$ stands for the set of measurable functions on  $I$ with values in $X$ such that $t\mapsto \|f(t)\|_X$ belongs to $L^p(I)$.

In the case where $I=[0,T],$  the space $L^p([0,T];X)$ (resp. $\cC([0,T];X)$ will also be denoted by $L_T^p(X)$ (resp. $\cC_T(X)$). Finally, if $I=\R^+$ we shall alternately use the notation $L^p(X).$

One of the main advantages of using the Littlewood-Paley decomposition is that it enables us to work with spectrally localized (hence smooth) functions rather than with rough objects. As a consequence, we naturally obtain bounds for each dyadic block in spaces of type $L^\rho_T(L^p).$  Going from those type of bounds to estimates in  $L^\rho_T(\dot B^s_{p,r})$ requires to perform a summation in $\ell^r(\Z).$ When doing so however, we \emph{do not} bound the $L^\rho_T(\dot B^s_{p,r})$ norm for the time integration has been performed \emph{before} the $\ell^r$ summation.
This leads to the following notation (after J.-Y. Chemin and N. Lerner in \cite{CL}):

\begin{defi}\label{d:espacestilde}
For $T>0,$ $s\in\R$ and  $1\leq r,\rho\leq\infty,$
 we set
$$
\|u\|_{\tilde L_T^\rho(\dot B^s_{p,r})}:=
\bigl\Vert2^{js}\|\ddq u\|_{L_T^\rho(L^p)}\bigr\Vert_{\ell^r(\Z)}.
$$
\end{defi}
One can then define the space $\tilde L^\rho_T(\dot B^s_{p,r})$ as the set of  tempered distributions $u$ over $(0,T)\times \R^d$ such that $\lim_{q\leftarrow-\infty}\dot S_q u=0$ in $L^\rho([0,T];L^\infty(\R^d))$ and $\|u\|_{\tilde L_T^\rho(\dot B^s_{p,r})}<\infty.$ The letter $T$ is omitted for functions defined over $\R^+.$ We shall also adopt the notation
$$
\tilde \cC_T(\dot B^s_{p,r}):=\tilde L_T^\infty(\dot B^s_{p,r})\cap \cC([0,T];\dot B^s_{p,r})
\quad\hbox{and}\quad 
\tilde\cC_b(\dot B^s_{p,r}):=\tilde L^\infty(\dot B^s_{p,r})\cap \cC_b(\R^+;\dot B^s_{p,r}).
$$
The spaces $\tilde L^\rho_T(\dot B^s_{p,r})$ may be compared with the spaces  $L_T^\rho(\dot B^s_{p,r})$ through the Minkowski inequality: we have
$$
\|u\|_{\tilde L_T^\rho(\dot B^s_{p,r})}
\leq\|u\|_{L_T^\rho(\dot B^s_{p,r})}\ \text{ if }\ r\geq\rho\quad\hbox{and}\quad
\|u\|_{\tilde L_T^\rho(\dot B^s_{p,r})}\geq
\|u\|_{L_T^\rho(\dot B^s_{p,r})}\ \text{ if }\ r\leq\rho.
$$
The general principle is that all the properties of continuity for the product and  composition which are true in Besov spaces (see below) remain
true in the above  spaces. The time exponent just behaves according to H\"older's inequality. 
\medbreak
Let us now recall a few nonlinear estimates in Besov spaces. Formally, any product of two distributions $u$ and $v$ may be decomposed into 
\begin{equation}\label{eq:bony}
uv=T_uv+T_vu+R(u,v)
\end{equation}
with: 
$$
T_uv:=\sum_q\dot S_{q-1}u\ddq v,\quad
T_vu:=\sum_q \dot S_{q-1}v\ddq u\ \hbox{ and }\ 
R(u,v):=\sum_q\sum_{|q'-q|\leq1}\ddq u\,\dot\Delta_{q'}v.
$$
The above operator $T$ is called ``paraproduct'' whereas $R$ is called ``remainder''. The decomposition \eqref{eq:bony} has been introduced by J.-M. Bony in \cite{BJM}.

Let us now briefly recall the following estimates (we refer to \cite{Dbook} section 2.6, \cite{Dinv}, \cite{Has1} for more properties of continuity for the paraproduct and remainder operators, sometimes adapted to $\tilde L_T^\rho(\dot B^s_{p,r})$ spaces). These estimates will be of constant use in the paper.

\begin{prop}
 \sl{
There exists a constant $C>0$ such that for any $s\in \R$ and any $p,r\in [1, \infty[$, we have for any $(u,v)\in L^\infty\times \dot{B}_{p,r}^s$,
$$
\|\dot{T}_u v\|_{\dot{B}_{p,r}^s}\leq C^{1+|s|} \|u\|_{L^\infty} \|v\|_{\dot{B}_{p,r}^s}.
$$
Moreover, for any $s\in \R$, any $t<0$, for any $p,r_1,r_2 \in [1,\infty[$, we have, for any $(u,v)\in \dot{B}_{\infty,r_1}^t\times \dot{B}_{p,r_2}^s$,
$$
\|\dot{T}_u v\|_{\dot{B}_{p,r}^{s+t}}\leq \frac{C^{1+|s+t|}}{-t} \|u\|_{\dot{B}_{\infty,r_1}^t} \|v\|_{\dot{B}_{p,r_2}^s} \quad \mbox{with} \quad \frac{1}{r}\overset{def}{=}\min(1, \frac{1}{r_1}+ \frac{1}{r_2}).
$$
}
\end{prop}
And concerning the homogeneous remainder:
\begin{prop}
\sl{ There exists a constant $C>0$ which satisfies the following inequalities. Let $s_1,s_2 \in \R$ and $p_1, p_2, r_1, r_2 \in [1,\infty]$. Let us assume that $s_1+s_2<\fdp$ or $s_1+s_2\leq\fdp$ if $r=1$, and
$$
\frac{1}{p}\overset{def}{=}\frac{1}{p_1}+ \frac{1}{p_2}\leq 1 \quad \mbox{and} \quad \frac{1}{r}\overset{def}{=}\frac{1}{r_1}+ \frac{1}{r_2}\leq 1.
$$
If $s_1+s_2>0$ then we have, for any $(u,v)\in \dot{B}_{p_1,r_1}^{s_1}.\dot{B}_{p_2,r_2}^{s_2}$,
$$
\|\dot{R}(u,v)\|_{\dot{B}_{p,r}^{s_1+s_2}} \leq \frac{C^{|s_1+s_2|+1}}{s_1+s_2} \|u\|_{\dot{B}_{p_1,r_1}^{s_1}} \|v\|_{\dot{B}_{p_2,r_2}^{s_2}}.
$$
When $r=1$ and $s_1+s_2\geq 0$, we have, for any $(u,v)\in \dot{B}_{p_1,r_1}^{s_1}.\dot{B}_{p_2,r_2}^{s_2}$,
$$
\|\dot{R}(u,v)\|_{\dot{B}_{p,\infty}^{s_1+s_2}} \leq C^{|s_1+s_2|+1} \|u\|_{\dot{B}_{p_1,r_1}^{s_1}} \|v\|_{\dot{B}_{p_2,r_2}^{s_2}}.
$$}
\end{prop}
Using Besov injections together with the previous estimates, we get the following result (we refer to \cite{Dbook} or \cite{Dmarne}):
\begin{prop}
\sl{There exists a constant $C>0$ which satisfies the following inequalities. For any $s_1, s_2 \in \R$, any $1\leq p_1,p_2,p \leq \infty$ and any $1\leq r_1,r_2,r \leq \infty$ such that
$$
s_1+s_2>0, \quad \frac{1}{p}\leq \frac{1}{p_1}+ \frac{1}{p_2}\leq 1, \quad \mbox{and} \quad \frac{1}{r}\leq \frac{1}{r_1}+ \frac{1}{r_2}\leq 1,
$$
we have
$$
\|\dot{R}(u,v)\|_{\dot{B}_{p,r}^{\sigma_{1,2}}} \leq \frac{C^{|s_1+s_2|+1}}{s_1+s_2} \|u\|_{\dot{B}_{p_1,r_1}^{s_1}} \|v\|_{\dot{B}_{p_2,r_2}^{s_2}}
\quad \mbox{with} \quad \sigma_{1,2}\overset{def}{=}s_1+s_2-d(\frac{1}{p_1}+\frac{1}{p_2}-\frac{1}{p}),
$$
provided that $\sigma_{1,2}<d/p$, or $\sigma_{1,2}\leq d/p$ and $r=1$.}
\end{prop}
In this article we will frequently use the following estimates (given by the previous recalls) under the same assumptions there exists constant $C>0$ such that:
$$
\|\dot{T}_u v\|_{\dot{B}_{2,1}^s}\leq C \|u\|_{L^\infty} \|v\|_{\dot{B}_{2,1}^s}\leq C \|u\|_{\dot{B}_{2,1}^\fd} \|v\|_{\dot{B}_{2,1}^s},
$$
$$\|\dot{T}_u v\|_{\dot{B}_{2,1}^{s+t}}\leq C\|u\|_{\dot{B}_{\infty,\infty}^t} \|v\|_{\dot{B}_{2,1}^s} \leq C\|u\|_{\dot{B}_{2,1}^{t+\fd}} \|v\|_{\dot{B}_{2,1}^s} \quad (t<0),
$$
$$\|\dot{R}(u,v)\|_{\dot{B}_{2,1}^{s_1+s_2}} \leq C\|u\|_{\dot{B}_{\infty,\infty}^{s_1}} \|v\|_{\dot{B}_{2,1}^{s_2}} \leq C\|u\|_{\dot{B}_{2,1}^{s_1+\fd}} \|v\|_{\dot{B}_{2,1}^{s_2}} \quad (s_1+s_2>0),
$$
\begin{equation}
 \|\dot{R}(u,v)\|_{\dot{B}_{2,1}^{s_1+s_2-\fd}} \leq C\|u\|_{\dot{B}_{2,1}^{s_1}} \|v\|_{\dot{B}_{2,\infty}^{s_2}} \leq C\|u\|_{\dot{B}_{2,1}^{s_1}} \|v\|_{\dot{B}_{2,1}^{s_2}} \quad (s_1+s_2>0).
\label{estimbesov}
\end{equation}
Let us now turn to the composition estimates. We refer for example to \cite{Dbook} (Theorem $2.59$, corollary $2.63$)):
\begin{prop}
\sl{\begin{enumerate}
 \item Let $s>0$, $u\in \dot{B}_{2,1}^s\cap L^{\infty}$ and $F\in W_{loc}^{[s]+2, \infty}(\R^d)$ such that $F(0)=0$. Then $F(u)\in \dot{B}_{2,1}^s$ and there exists a function of one variable $C_0$ only depending on $s$, $d$ and $F$ such that
$$
\|F(u)\|_{\dot{B}_{2,1}^s}\leq C_0(\|u\|_{L^\infty})\|u\|_{\dot{B}_{2,1}^s}.
$$
\item If $u$ and $v\in\dot{B}_{2,1}^\fd$ and if $v-u\in \dot{B}_{2,1}^s$ for $s\in]-\fd, \fd]$ and $G\in W_{loc}^{[s]+3, \infty}(\R^d)$ such that $G'(0)=0$, then $G(v)-G(u)$ belongs to $\dot{B}_{2,1}^s$ and there exists a function of two variables $C$ only depending on $s$, $d$ and $G$ such that
$$
\|G(v)-G(u)\|_{\dot{B}_{2,1}^s}\leq C(\|u\|_{L^\infty}, \|v\|_{L^\infty})\left(\|u\|_{\dot{B}_{2,1}^\fd} +\|v\|_{\dot{B}_{2,1}^\fd}\right) \|v-u\|_{\dot{B}_{2,1}^s}.
$$
\end{enumerate}}
\label{estimcompo}
\end{prop}
\begin{rem}
 \sl{As $2.$ is proved using $1.$, together with the fact that $G(v)-G(u)=(v-u).\int_0^1 G'(u+\tau(v-u))d\tau$ and the previous estimates on the paraproduct and remainder, we can prove that, more generally,
$$
\|G(v)-G(u)\|_{\dot{B}_{2,1}^s}\leq C(\|u\|_{L^\infty}, \|v\|_{L^\infty})\left(|G'(0)| +\|u\|_{\dot{B}_{2,1}^\fd} +\|v\|_{\dot{B}_{2,1}^\fd}\right) \|v-u\|_{\dot{B}_{2,1}^s}.
$$
}
\end{rem}

\subsubsection{Hybrid Besov spaces}

As announced before we will get a different regularity for the density fluctuation in some low and high frequencies, separated by a frequency threshold. This will define the hybrid Besov spaces. Let us begin with the spaces that are introduced by R. Danchin in \cite{Dinv} or \cite{Dbook}:

\begin{defi}
 \sl{
For $\alpha>0$, $r\in [0, \infty]$ and $s\in \R$ we denote
$$
\|u\|_{\tilde{B}_\alpha^{s,r}} \overset{def}{=} \Sum_{l\in \Z} 2^{ls} \max(\alpha, 2^{-l})^{1-\frac{2}{r}}\|\ddl u\|_{L^2}
$$}
\end{defi}
In the present paper, we will only use these norms with $r\in \{1,\infty\}$:
$$
\|u\|_{\tilde{B}_\alpha^{s,\infty}}= \Sum_{l\leq \log_2(\frac{1}{\alpha})} 2^{l(s-1)} \|\ddl u\|_{L^2}+ \Sum_{l> \log_2(\frac{1}{\alpha})} \alpha 2^{ls} \|\ddl u\|_{L^2},
$$
and
$$
\|u\|_{\tilde{B}_\alpha^{s,1}}= \Sum_{l\leq \log_2(\frac{1}{\alpha})} 2^{l(s+1)} \|\ddl u\|_{L^2}+ \Sum_{l> \log_2(\frac{1}{\alpha})} \frac{1}{\alpha} 2^{ls} \|\ddl u\|_{L^2},
$$
\begin{rem}
\sl{As stated in \cite{Dbook} we have the equivalence
$$
\frac{1}{2} \left( \|u\|_{\dot{B}_{2,1}^{s-1}}+ \alpha \|u\|_{\dot{B}_{2,1}^s}\right)\leq \|u\|_{\tilde{B}_\alpha^{s,\infty}} \leq \|u\|_{\dot{B}_{2,1}^{s-1}}+ \alpha \|u\|_{\dot{B}_{2,1}^s}.
$$
And if we denote $\|u\|_{\dot{B}_{2,1}^{s-1}\cap \dot{B}_{2,1}^s}=\|u\|_{\dot{B}_{2,1}^{s-1}}+ \|u\|_{\dot{B}_{2,1}^s}$, then
$$
\frac{1}{2}\min(1,\alpha) \|u\|_{\dot{B}_{2,1}^{s-1}\cap \dot{B}_{2,1}^s}\leq \|u\|_{\tilde{B}_\alpha^{s,\infty}} \leq \max(1,\alpha) \|u\|_{\dot{B}_{2,1}^{s-1}\cap \dot{B}_{2,1}^s}
$$} 
\end{rem}
In this article we will mainly use the following hybrid Besov norm:
\begin{defi}
\sl{If $l_\ee$ is a frequency threshold, and $s,t\in \R$ we define the following hybrid norms:
\begin{equation}
\|u\|_{\dot{B}_{\ee}^{s,t}} \overset{def}{=} \Sum_{l\leq l_\ee} 2^{ls} \|\ddl u\|_{L^2}+ \Sum_{l> l_\ee} \frac{1}{\ee^2} 2^{lt} \|\ddl u\|_{L^2}.
 \label{normhybride}
\end{equation}}
\end{defi}
\begin{rem}
\sl{In the following we will take $l_\ee=[\frac{1}{2}\log_2(\frac{\gamma}{C_0\ee^2})-1]$ as defined in (\ref{freqseuil}). It goes to infinity as $\ee$ goes to zero.}
\end{rem}
\begin{rem}
 \sl{If $\ee$ is small enough, we have $\|u\|_{\tilde{B}_1^{s,1}} \leq \|u\|_{\tilde{B}_\ee^{s,1}} \leq \|u\|_{\dot{B}_{\ee}^{s+1,s}}$}
\label{hybridcomparaison}
\end{rem}

\subsection{The Korteweg System}

Unlike the Navier-Stokes or $(NSR)$ systems, here the density, as the velocity, is regularized for every frequency (in fact there is a frequency threshold in the Fourier modes, but in both cases the density is parabolically regularized).  In the sequel we choose $\overline{\rho}=1$ and we will denote $q_0=\rho_0-1$ and $q=\rho-1$ the density fluctuations. We will recall here some results from \cite{DD} about the strong solutions of the Korteweg system (for more simplicity, we do not mention here any exterior forcing term):

\begin{thm}
\sl{ Assume that $P'(1)>0$, $\min(\mu, 2 \mu+\lambda)>0$, that the initial density fluctuation $q_0$ belongs to $\dot{B}_{2,1}^{\fd-1}\cap \dot{B}_{2,1}^{\fd}$, and that the initial velocity $u_0$ is in $(\dot{B}_{2,1}^{\fd-1})^d$. Then there exist constants $\eta_K>0$ and $C>0$ depending on $\kappa$, $\mu$, $\nu$, $P'(1)$ and $d$ such that if:
$$
\|q_0\|_{\dot{B}_{2,1}^{\fd-1} \cap \dot{B}_{2,1}^{\fd}} +\|u_0\|_{\dot{B}_{2,1}^{\fd-1}}\leq \eta_K
$$
then system $(NSK)$ has a unique global solution $(\rho, u)$ such that the density fluctuation and the velocity satisfy:
$$
\begin{cases}
q \in \cC(\R_+, \dot{B}_{2,1}^{\fd-1} \cap \dot{B}_{2,1}^{\fd}) \cap L^1(\R_+, \dot{B}_{2,1}^{\fd+1} \cap \dot{B}_{2,1}^{\fd+2}),\\
u\in \cC(\R_+, \dot{B}_{2,1}^{\fd-1})^d \cap  L^1(\R_+, \dot{B}_{2,1}^{\fd+1})^d.
\end{cases}
$$
Moreover the norm of $(q,u)$ in this space is estimated by the initial norm $C (\|q_0\|_{\dot{B}_{2,1}^{\fd-1} \cap \dot{B}_{2,1}^{\fd}} +\|u_0\|_{\dot{B}_{2,1}^{\fd-1}})$.}
\end{thm}
Further in this article R. Danchin and B. Desjardins provide a Fourier study of the linearized system and observe different behaviours wether the quantity $\nu^2- 4\kappa$ is positive, negative of zero. In any case they obtain parabolic regularization.

\subsection{Notations and main results}


\begin{defi}
\sl{The space $E_{\ee}^s$ is the set of functions $(q,u)$ in
$$
\left(\cC_b(\R_+, \dot{B}_{2,1}^{s-1}\cap \dot{B}_{2,1}^s)\cap L^1(\R_+, \dot{B}_\ee^{s+1,s}\cap \dot{B}_\ee^{s+2,s})\right) \times
\left(\cC_b(\R_+, \dot{B}_{2,1}^{s-1})\cap L^1(\R_+, \dot{B}_{2,1}^{s+1})\right)^d
$$
endowed with the norm
\begin{multline}
\|(q,u)\|_{E_{\ee}^s} \overset{def}{=} \|u\|_{L^{\infty} \dot{B}_{2,1}^{s-1}}+ \|q\|_{L^{\infty} \dot{B}_{2,1}^{s-1}}+ \|q\|_{L^{\infty} \dot{B}_{2,1}^{s}}\\
+\|u\|_{L^1 \dot{B}_{2,1}^{s+1}}+ \|q\|_{L^1 \dot{B}_{\ee}^{s+1,s}}+ \|q\|_{L^1 \dot{B}_{\ee}^{s+2,s}}
\end{multline}}
\end{defi}
We first prove global well-posedness for system $(NSR_\ee)$ by following similar ideas as in \cite{Has1} and uniform estimates with respect to $\ee$:
\begin{thm}
\sl{Let $\ee>0$ and assume that $\min(\mu,2\mu+\lambda)>0$. There exist two positive constants $\eta_R$ and $C$ only depending on $d$, $\kappa$, $\mu$, $\lambda$ and $P'(1)$ such that if $q_0\in \dot{B}_{2,1}^{\fd-1}\cap \dot{B}_{2,1}^{\fd}$, $u_0 \in \dot{B}_{2,1}^{\fd-1}$ and
$$
\|q_0\|_{\dot{B}_{2,1}^{\fd-1} \cap \dot{B}_{2,1}^{\fd}} +\|u_0\|_{\dot{B}_{2,1}^{\fd-1}}\leq \eta_R
$$
then system $(NSR_\ee)$ has a unique global solution $(\rho, u)$ with $(q, u)\in E_{\ee}^{\fd}$ such that:
$$
\|(q,u)\|_{E_{\ee}^\fd} \leq C (\|q_0\|_{\dot{B}_{2,1}^{\fd-1} \cap \dot{B}_{2,1}^{\fd}} +\|u_0\|_{\dot{B}_{2,1}^{\fd-1}}).
$$}
\label{thexist}
\end{thm}
\begin{rem}
 \sl{Note that in the low frequencies, the parabolic regularization for $q$ is the same as for the Korteweg system, that is the low frequencies of $q$ are in $\dot{B}_{2,1}^{\fd+1} \cap \dot{B}_{2,1}^{\fd+2}$}.
\end{rem}

The main result in this article is the following: when the initial data are small enough (so that we have global solutions for $(NSK)$ and $(NSR_\ee)$) the solution of $(NSR_\ee)$ goes to the solution of $(NSK)$ when $\ee$ goes to zero.

\begin{thm}
\sl{Assume that $\min(\mu,2\mu+\lambda)>0$, $P'(1)>0$ and that $q_0\in \dot{B}_{2,1}^{\fd-1}\cap \dot{B}_{2,1}^{\fd}$, $u_0 \in \dot{B}_{2,1}^{\fd-1}$. There exists $0<\eta\leq \min(\eta_K, \eta_R)$ such that if
$$
\|q_0\|_{\dot{B}_{2,1}^{\fd-1} \cap \dot{B}_{2,1}^{\fd}} +\|u_0\|_{\dot{B}_{2,1}^{\fd-1}}\leq \eta,
$$
then systems $(NSK)$ and $(NSR_\ee)$ both have global solutions, and with the same notations as before, there exists a constant $C=C(\eta, \kappa, \overline{\rho}, P'(1))>0$ such that for all $\alpha\in ]0, 1[$ (if $d=2$) or $\alpha\in ]0,1]$ (if $d\geq 3$),
\begin{multline}
 \|\ue-u\|_{\tilde{L}_t^{\infty} \dot{B}_{2,1}^{\fd-\alpha-1}}+ \|\qe-q\|_{\tilde{L}_t^{\infty} \dot{B}_{2,1}^{\fd-\alpha-1}}+ \|\qe-q\|_{\tilde{L}_t^{\infty} \dot{B}_{2,1}^{\fd-\alpha}}\\
+\|\ue-u\|_{\tilde{L}_t^1 \dot{B}_{2,1}^{\fd-\alpha+1}}+ \|\qe-q\|_{\tilde{L}_t^1 \dot{B}_{\ee}^{\fd-\alpha+1,\fd-\alpha}}+ \|\qe-q\|_{\tilde{L}_t^1 \dot{B}_{\ee}^{\fd-\alpha+2,\fd-\alpha}} \leq C \ee^{\alpha}.
\end{multline}}
\label{thcv}
\end{thm}

\begin{rem}
 \sl{The same results hold for any function $\phi\in \cS(\R^d)$ such that $\forall \xi\in \R^d$, $\hat{\phi}(\xi)=g(|\xi|^2)$ with:
\begin{itemize}
 \item Function $g:\R_+\rightarrow \R$ takes its values in $[0,1]$, with $g(0)=1$,
 \item Function $h:x\mapsto \frac{1-g(x)}{x}$ is decreasing with $\lim_0 h=1$ and $\lim_\infty h=0$,
 \item Function $k:x\mapsto 1-g(x)$ is increasing with $\lim_0 k=0$ and $\lim_\infty k=1$,
 \item For all $1<\beta<2$, there exists $C_\beta>0$ such that for all $x\geq 0$,
$$
0\leq \frac{g(x)-1+x}{x^\beta}\leq C_\beta.
$$
\end{itemize}
}
\end{rem}

The structure of this article is the following: in the second section we will obtain a priori estimates on the linear system with convection terms. The third section is devoted to the existence and uniqueness of solutions for the non-local model, and in the last section we will obtain the convergence result of theorem \ref{thcv}. In the appendix, one will find the proof of some estimates involving the special hybrid Besov norms introduced in this paper.

\section{A priori estimates}

In this section we focus on the following linear system ($\ee>0$ is fixed and for more simplicity we write $(q,u)$ instead of $(\qe,\ue)$):
$$
\begin{cases}
\begin{aligned}
&\d_t q+ v.\nabla q+ \div u= F,\\
&\d_t u+ v.\nabla u -\cA u+ p\nabla q-\frac{k}{\ee^2} \nabla(\phi_\ee*q-q)= G.\\
\end{aligned}
\end{cases}
\leqno{(LR_\ee)}
$$
With
$$\cA u= \mu \Delta u+ (\lambda+\mu)\nabla \div u.$$
This section is devoted to the proof of the following a priori estimates:
\begin{prop} Let $\ee>0$, $s\in \R$, $I=[0,T[$ or $[0, +\infty[$ and $v\in L^1(I,\dot{B}_{2,1}^{\fd+1}) \cap L^2 (I,\dot{B}_{2,1}^{\fd})$. Assume that $(q,u)$ is a solution of System $(LR_\ee)$ defined on $I$. There exists a constant $C>0$ depending on $d$, $s$, $\mu$, $\nu$, $p$, $k$, $c_0$ and $C_0$ such that for all $t\in I$,
\begin{multline}
 \|u\|_{\tilde{L}_t^{\infty} \dot{B}_{2,1}^{s-1}}+ \|q\|_{\tilde{L}_t^{\infty} \dot{B}_{2,1}^{s-1}}+ \|q\|_{\tilde{L}_t^{\infty} \dot{B}_{2,1}^{s}}+ \|u\|_{\tilde{L}_t^1 \dot{B}_{2,1}^{s+1}}+ \|q\|_{\tilde{L}_t^1 \dot{B}_{\ee}^{s+1,s}}+ \|q\|_{\tilde{L}_t^1 \dot{B}_{\ee}^{s+2,s}}\\
\leq C e^{C\int_0^t (\|\nabla v(\tau)\|_{\dot{B}_{2,1}^\fd}+ \|v(\tau)\|_{\dot{B}_{2,1}^\fd}^2)d\tau} \Big(\|u_0\|_{\dot{B}_{2,1}^{s-1}}+ \|q_0\|_{\dot{B}_{2,1}^{s-1}} + \|q_0\|_{\dot{B}_{2,1}^{s}}\\
+ \|F\|_{\tilde{L}_t^1 \dot{B}_{2,1}^{s-1}}+ \|F\|_{\tilde{L}_t^1 \dot{B}_{2,1}^{s}}+ \|G\|_{\tilde{L}_t^1 \dot{B}_{2,1}^{s-1}}\Big).
\label{estimapriori}
\end{multline}
\label{apriori}
\end{prop}
The proof of the proposition will be close to those in \cite{Dinv}, \cite{Dbook}, \cite{Dmach} or \cite{Has1}: lead by the behaviour of the linearized system and using symmetrizers we get estimates on dyadic blocks and obtain the expected estimates. The difference here is that we need to carefully estimate penalized terms and we obtain, for the density, a frequency threshold depending on $\ee$.

We choose here not to separate the compressible and incompressible parts, and we will localize in frenquency: using the Littlewood-Paley decomposition, for all $l\in \Z$ we define $\ql= \ddl q$ and $\ul= \ddl u$. We obtain the following system:
$$
\begin{cases}
\begin{aligned}
&\d_t \ql+ v.\nabla \ql+ \div \ul= \Fl+\Rl,\\
&\d_t \ul+ v.\nabla \ul -\cA \ul+ p\nabla \ql-\frac{k}{\ee^2} \nabla(\phi_\ee*\ql-\ql)= \Gl+\Rlp,\\
\end{aligned}
\end{cases}
\leqno{(LLR_\ee)}
$$
where $\Rl= [v.\nabla, \ddl] q$ and $\Rlp= [v.\nabla, \ddl] u$.

If $\alpha>0$ is fixed and small enough (for example we can choose $\alpha=\min(\mu, \lambda+2\mu)/4$), we will introduce in the following:
$$
\hl^2=\|\ul\|_{L^2}^2+ p \|\ql\|_{L^2}^2+ \alpha \left(\nu \|\nabla \ql\|_{L^2}^2+ 2(\ul|\nabla \ql)_{L^2}\right) +\frac{k}{\ee^2}(\ql|\ql- \phi_\ee* \ql)_{L^2}.
$$

\subsection{Preliminary result}
In order to prove the estimate we will show as a first step the following lemma:

\begin{lem}
Under the previous notations, there exist some $m>0$ and $C>0$ such that for all $l\in \Z$:
\begin{equation}
\begin{aligned}
&\frac{1}{2} \frac{d}{dt} \hl^2+ m\big(2^{2l} \|\ul\|_{L^2}^2 + \|\nabla \ql\|_{L^2}^2 +\frac{1}{\ee^2} (\nabla\ql|\nabla\ql- \phi_\ee* \nabla\ql)_{L^2}\big)\\
&\hspace{1cm}\leq C (\|\nabla v\|_{L^\infty}+ \|v\|_{L^\infty}^2)\hl^2 + C\big((1+2^l)(\|\Fl\|_{L^2}+ \|\Rl\|_{L^2})\\
&\hspace{8cm}+ \|\Gl\|_{L^2}+ \|\Rlp\|_{L^2}\big) \hl.
\end{aligned}
\end{equation}
\label{energieZ}
\end{lem}
\textbf{Proof: } If we compute the innerproduct in $L^2$ of the first equation from system $(LLR_\ee)$ by $\ql$, and the innerproduct of the second equation by $\ul$ we obtain that for all $l\in \Z$: 
\begin{equation}
  \Frac{1}{2} \Frac{d}{dt} \|\ql\|_{L^2}^2+ (\div \ul|\ql)_{L^2}=-(v.\nabla \ql|\ql)_{L^2}+ (\Fl+\Rl|\ql)_{L^2},
\end{equation}
and
\begin{multline}
 \Frac{1}{2} \Frac{d}{dt} \|\ul\|_{L^2}^2-(\cA \ul|\ul)_{L^2}+p(\nabla \ql|\ul)_{L^2} -\frac{k}{\ee^2} (\phi_\ee* \nabla\ql-\nabla\ql|\ul)_{L^2}\\
= -(v.\nabla \ul|\ul)_{L^2}+ (\Gl+\Rlp|\ul)_{L^2}.
\end{multline}
Using an integration by parts allows us to write that $(\div \ul|\ql)_{L^2}= -(\nabla \ql|\ul)_{L^2}$ and then to combine the previous estimates in order to eliminate these terms. We need to get rid of terms like this because they cannot be absorbed by the left-hand side and after the use of a Gronwall type estimate they introduce multiplicative constant terms $e^{Ct}$ (which are problematic as we look for global time estimates). It is less easy to get rid of the term $\Frac{k}{\ee^2} (\phi_\ee* \nabla\ql-\nabla\ql|\ul)_{L^2}$ and we will explain later how to do it.\\

Integrations by parts also provide that there exists a constant $C>0$ such that:
$$
\begin{cases}
 |(v.\nabla \ql|\ql)_{L^2}|\leq C \|\nabla v\|_{L^\infty} \|\ql\|_{L^2}^2,\\
|(v.\nabla \ul|\ul)_{L^2}|\leq C \|\nabla v\|_{L^\infty} \|\ul\|_{L^2}^2,
\end{cases}
$$
as well as the fact that
$$-(\cA \ul|\ul)_{L^2}= \mu \|\nabla \ul\|_{L^2}^2+(\lambda+ \mu)\|\div \ul\|_{L^2}^2$$
which leads to
$$-(\cA \ul|\ul)_{L^2}\geq \underline{\nu} \|\nabla \ul\|^2, \quad \mbox{where} \quad \nu=\lambda+2 \mu \quad \mbox{and} \quad \underline{\nu}=\min(\mu, \nu).$$
(if $\lambda+\mu\geq 0$ it is immediate, else we use that $\|\div \ul\|_{L^2}^2\leq \|\nabla \ul\|^2$.)
Combining these estimates we obtain that:
\begin{multline}
 \Frac{1}{2} \Frac{d}{dt} (\|\ul\|_{L^2}^2+ p\|\ql\|_{L^2}^2) +\underline{\nu} \|\nabla \ul\|_{L^2}^2 -\frac{k}{\ee^2} (\phi_\ee* \nabla\ql-\nabla\ql|\ul)_{L^2}\\
\leq C\|\nabla v\|_{L^\infty} (\|\ul\|_{L^2}^2+ p\|\ql\|_{L^2}^2) +p(\|\Fl\|_{L^2}+\|\Rl\|_{L^2})\|\ql\|_{L^2} \\
+(\|\Gl\|_{L^2}+\|\Rlp\|_{L^2})\|\ul\|_{L^2}.
\label{equq}
\end{multline}
As we want to obtain regularization on the density, we need at least a term such as $\|\ql\|_{L^2}^2$ to appear in the left-hand side. If we compute the inner product of the equation on $\ul$ by $\nabla \ql$ we get:
\begin{multline}
 (\d_t \ul| \nabla \ql)_{L^2}+(v.\nabla \ul|\nabla \ql)_{L^2}-(\cA \ul| \nabla \ql)_{L^2}+p\|\nabla \ql\|_{L^2}^2+\frac{k}{\ee^2}(\nabla \ql- \phi_\ee*\nabla \ql|\nabla \ql)_{L^2}\\
=(\Gl+ \Rlp|\nabla \ql)_{L^2}.
\label{en1}
\end{multline}

In order to get a full derivative, we need to estimate the term $(\ul|\d_t \nabla \ql)$ so we will write the equation on $\nabla \ql$ and its innerproduct on $\ul$:
\begin{equation}
 \d_t \nabla \ql+ \nabla(v.\nabla \ql)+ \nabla \div \ul= \nabla \Fl+ \nabla \Rl.
\end{equation}
Before adding (\ref{en1}) to
\begin{equation}
 (\d_t \nabla \ql|\ul)_{L^2}+ (\nabla(v.\nabla \ql)|\ul)_{L^2}+ (\nabla \div \ul|\ul)_{L^2}= (\nabla \Fl+ \nabla \Rl|\ul)_{L^2},
\end{equation}
we remark that a simple computation shows there exists a constant $C>0$ such that:
$$
|(v.\nabla \ul|\nabla \ql)_{L^2}+(\nabla(v.\nabla \ql)|\ul)_{L^2}|\leq C \|\nabla v\|_{L^\infty} \|\nabla \ql\|_{L^2}\|\ul\|_{L^2}.
$$
Moreover, integrations by part also provide the fact that
$$
(\cA \ul| \nabla \ql)_{L^2}= \nu (\ul| \nabla \Delta \ql)_{L^2}.
$$
Gathering these estimates, we obtain that:
\begin{multline}
 \frac{d}{dt} (\ul| \nabla \ql)_{L^2}-\nu (\ul| \nabla \Delta \ql)_{L^2} -\|\div \ul\|_{L^2}^2+ p\|\nabla \ql\|_{L^2}^2 +\frac{k}{\ee^2}(\nabla \ql- \phi_\ee*\nabla \ql|\nabla \ql)_{L^2}\\
\leq C \|\nabla v\|_{L^\infty} \|\nabla \ql\|_{L^2}\|\ul\|_{L^2}+ (\|\Gl\|_{L^2}+ \|\Rlp\|_{L^2})\|\nabla \ql\|_{L^2}+ 2^l (\|\Fl\|_{L^2}+ \|\Rl\|_{L^2})\|\ul\|_{L^2}.
\label{eqmixte}
\end{multline}
\begin{rem}
\sl{Note that $\Frac{k}{\ee^2}(\nabla \ql- \phi_\ee*\nabla \ql|\nabla \ql)_{L^2}$ is nonnegative.}
\label{absorb}
\end{rem}
In the previous equation, the term $(\ul| \nabla \Delta \ql)_{L^2}$ is annoying because it introduces much more derivatives that we will be able to handle in the following, so we need to neutralize it. A simple way to do that is to compute the innerproduct of the equation of $\nabla \ql$ by $\nabla \ql$. Using that:
$$
(\nabla \div \ul| \nabla \ql)_{L^2}=(\ul| \nabla \Delta \ql)_{L^2},
$$
and that
$$
|(\nabla(v.\nabla \ql)| \nabla \ql)_{L^2}|\leq C \|\nabla v\|_{L^\infty} \|\nabla \ql\|_{L^2}^2, 
$$
we get the following estimate:
\begin{multline}
  \Frac{1}{2} \Frac{d}{dt} \|\nabla \ql\|_{L^2}^2+ (\ul| \nabla \Delta \ql)_{L^2} \leq C \|\nabla v\|_{L^\infty} \|\nabla \ql\|_{L^2}^2+ 2^l(\|\Fl\|_{L^2}+ \|\Rl\|_{L^2})\|\nabla\ql\|_{L^2},
\label{eqdq}
\end{multline}
so that, adding $\nu$ (\ref{eqdq}) to (\ref{eqmixte}) will neutralize the previous annoying term.
\begin{rem}
 \sl{Due to $-\|\div \ul\|_{L^2}^2$ appearing in (\ref{eqmixte}) we have to be careful and look at (\ref{equq})+$\alpha$ ($\nu$ (\ref{eqdq})+(\ref{eqmixte})) with $\alpha>0$ small enough ($\alpha< \underline{\nu}$) so that $-\alpha \|\div \ul\|_{L^2}^2$ can be absorbed by $\underline{\nu}\|\nabla \ul\|_{L^2}^2$.}
\end{rem}
As announced we have now to manage to get rid of $\Frac{k}{\ee^2} (\phi_\ee* \nabla\ql-\nabla\ql|\ul)_{L^2}$ which is introduced by estimate (\ref{equq}). A rough estimate will eventually provide some constant term $\ee^{-2}$ which will ruin any asymptotic approach as $\ee$ goes to zero. The easiest way to take care of it is to consider the equation satisfied by the nonnegative term $\Frac{k}{\ee^2}(\ql|\ql-\phi_{\ee}*\ql)_{L^2}$. Computing the innerproduct of the equation on $\ql$ by $\ql-\phi_{\ee}*\ql$ gives:
\begin{multline}
 \Frac{1}{2} \Frac{d}{dt}(\ql|\ql-\phi_{\ee}*\ql)_{L^2}-(\ul|\nabla\ql-\phi_\ee* \nabla\ql)_{L^2}\\
=-(v.\nabla \ql|\ql-\phi_{\ee}*\ql)_{L^2}+(\Fl+ \Rl|\ql-\phi_{\ee}*\ql)_{L^2}.
\label{correct}
\end{multline}
We need to be careful when estimating the right-hand side as this estimate will be multiplied by $k/\ee^2$ and added to the others. Using Plancherel (up to a multiplicative constant) we introduce the following notation:
$$
(\ql| \ql-\phi_{\ee}*\ql)_{L^2}=\g |\hat{\ql}(\xi)|^2 (1-e^{-\ee^2 |\xi|^2})d\xi \overset{\mbox{def}}{=} \|\sqrt{1-\phi_\ee}*\ql\|_{L^2}^2.
$$
When estimating the right-hand side, we are forced to let the block $\ql-\phi_{\ee}*\ql$ in one part, as our only possibility is to absorb some terms using remark \ref{absorb}. As the following function:
\begin{equation}
\begin{aligned}
 {f:\quad \R_+^*} & {\to\R}\\
{x} & {\mapsto \frac{1-e^{-x}}{x}},
\end{aligned}
\label{foncaux}
\end{equation}
is decreasing from $1$ to $0$, and as up to loosing some derivatives we want to get uniform estimates in $\ee$, we can write that for all $\xi \in \R^d$:
$$
1-e^{-\ee^2 |\xi|^2}\leq \ee^2 |\xi|^2.
$$
So, using Plancherel,
\begin{multline}
 |(v.\nabla \ql|\ql-\phi_{\ee}*\ql)_{L^2}|=C |\g \hat{v.\nabla \ql}.\sqrt{1-e^{-\ee^2 |\xi|^2}}.\sqrt{1-e^{-\ee^2 |\xi|^2}}.\overline{\hat{\ql}(\xi)} d\xi|\\
 \leq C \g |\hat{v.\nabla \ql}|.\ee|\xi|.\sqrt{1-e^{-\ee^2 |\xi|^2}}.|\overline{\hat{\ql}(\xi)}| d\xi.
\end{multline}
Thanks to the fact that we are localized in frequency, and using the Cauchy-Schwarz estimate, we obtain that
\begin{multline}
  |(v.\nabla \ql|\ql-\phi_{\ee}*\ql)_{L^2}|\leq C \ee 2^l \|v.\nabla \ql\|_{L^2} \|\sqrt{1-\phi_\ee}*\ql\|_{L^2}^2\\
\leq C \ee 2^l \|v\|_{L^{\infty}} \|\nabla \ql\|_{L^2} \sqrt{(\ql| \ql-\phi_{\ee}*\ql)_{L^2}}.
\label{transp}
\end{multline}
Using the same arguments allows us to write:
\begin{multline}
|(\Fl+ \Rl|\ql-\phi_{\ee}*\ql)_{L^2}|=C|\g \hat{\Fl+ \Rl}.\sqrt{1-e^{-\ee^2 |\xi|^2}}.\sqrt{1-e^{-\ee^2 |\xi|^2}}.\overline{\hat{\ql}(\xi)}d\xi|\\
\leq C \ee 2^l (\|\Fl\|_{L^2}+\|\Rl\|_{L^2})\sqrt{(\ql| \ql-\phi_{\ee}*\ql)_{L^2}}.
\label{forcext}
\end{multline}
Plugging (\ref{transp}) and (\ref{forcext}) into (\ref{correct}), we obtain that:
\begin{multline}
 \Frac{1}{2} \Frac{d}{dt}(\ql|\ql-\phi_{\ee}*\ql)_{L^2}-(\ul|\nabla\ql-\phi_\ee* \nabla\ql)_{L^2}
 \leq C \ee^2 \|v\|_{L^{\infty}} \|\nabla \ql\|_{L^2} \sqrt{\frac{(\nabla\ql|\nabla\ql-\phi_{\ee}*\nabla\ql)_{L^2}}{\ee^2}}\\
+C \ee^2 2^l (\|\Fl\|_{L^2}+\|\Rl\|_{L^2})\sqrt{\frac{(\ql| \ql-\phi_{\ee}*\ql)_{L^2}}{\ee^2}}.
\label{correct2}
\end{multline}
Now, when we compute (\ref{equq})$+\alpha$ ($\nu$ (\ref{eqdq})$+$(\ref{eqmixte}))$+\frac{k}{\ee^2}$(\ref{correct2}), we will be able to neutralize or absorb any annoying term. For that we will introduce as previously announced:
$$
\hl^2=\|\ul\|_{L^2}^2+ p \|\ql\|_{L^2}^2+ \alpha \left(\nu \|\nabla \ql\|_{L^2}^2+ 2(\ul|\nabla \ql)_{L^2}\right) +\frac{k}{\ee^2}(\ql|\ql- \phi_\ee* \ql)_{L^2}.
$$
Using the Cauchy-Schwarz estimate, we can write that
$$
|(\ul|\nabla \ql)|_{L^2}\leq \|\ul\|_{L^2}.\|\nabla \ql\|_{L^2}\leq \frac{\nu}{4}\|\nabla \ql\|_{L^2}^2 +\frac{1}{\nu}\|\ul\|_{L^2}^2,
$$
so
\begin{multline}
 (1-\frac{2\alpha}{\nu})\|\ul\|_{L^2}^2+ p\|\ql\|_{L^2}^2 + \frac{\alpha \nu}{2} \|\nabla\ql\|_{L^2}^2+ \frac{k}{\ee^2}(\ql|\ql-\phi_\ee* \ql)_{L^2} \leq \hl^2\\
\leq (1+\frac{2\alpha}{\nu})\|\ul\|_{L^2}^2+ p\|\ql\|_{L^2}^2 + \frac{3\alpha \nu}{2} \|\nabla\ql\|_{L^2}^2+ \frac{k}{\ee^2}(\ql|\ql-\phi_\ee* \ql)_{L^2}.
\label{hl}
\end{multline}
The last term is nonnegative and lesser than $\|\nabla\ql\|_{L^2}$. If we want $\hl^2$ to be equivalent to $\|\ul\|_{L^2}^2+ \|\ql\|_{L^2}^2 +\|\nabla\ql\|_{L^2}^2$ we need that $\alpha< \nu/2$. Remember that we already have a condition on $\alpha$: for $-\alpha \|\div \ul\|_{L^2}^2$ to be absorbed by the left-hand side, we need $\alpha< \underline{\nu}$. So we need $\alpha<\min(\mu, \frac{\nu}{2})$ and for example we will fix here once and for all $\alpha=\frac{\underline{\nu}}{4}$. After computing (\ref{equq})$+\alpha$ ($\nu$ (\ref{eqdq})$+$(\ref{eqmixte}))$+\frac{k}{\ee^2}$(\ref{correct2}) we get:
\begin{multline}
 \Frac{1}{2} \Frac{d}{dt} \hl^2+ (\underline{\nu}-\alpha)\|\nabla \ul\|_{L^2}^2+ \alpha p \|\nabla \ql\|_{L^2}^2+ \frac{\alpha k}{\ee^2}(\nabla \ql- \phi_\ee*\nabla \ql|\nabla \ql)_{L^2}\\
\leq C \|\nabla v\|_{L^{\infty}} \hl^2 +C\left((1+2^l)(\|\Fl\|_{L^2}+ \|\Rl\|_{L^2})+ \|\Gl\|_{L^2}+ \|\Rlp\|_{L^2}\right) \hl\\
+C k 2^l (\|\Fl\|_{L^2}+\|\Rl\|_{L^2})\sqrt{\frac{(\ql| \ql-\phi_{\ee}*\ql)_{L^2}}{\ee^2}}\\
+C k \|v\|_{L^{\infty}} \|\nabla \ql\|_{L^2} \sqrt{\frac{(\nabla\ql|\nabla\ql-\phi_{\ee}*\nabla\ql)_{L^2}}{\ee^2}}.
\end{multline}
Thanks to the definition of $\hl$, we have:
$$
k 2^l (\|\Fl\|_{L^2}+\|\Rl\|_{L^2})\sqrt{\frac{(\ql| \ql-\phi_{\ee}*\ql)_{L^2}}{\ee^2}}\leq \sqrt{k} 2^l (\|\Fl\|_{L^2}+\|\Rl\|_{L^2}) \hl,
$$
and using the classical estimate $ab\leq (a^2+ b^2)/2$, we obtain for the last term:
\begin{multline}
 C k \|v\|_{L^{\infty}} \|\nabla \ql\|_{L^2} \sqrt{\frac{(\nabla\ql|\nabla\ql-\phi_{\ee}*\nabla\ql)_{L^2}}{\ee^2}}\\
\leq \frac{\alpha k}{2\ee^2}(\nabla\ql|\nabla\ql-\phi_{\ee}*\nabla\ql)_{L^2}+ \frac{k C^2}{2 \alpha} \|v\|_{L^{\infty}}^2 \|\nabla \ql\|_{L^2}^2.
\end{multline}
So, if we denote by $m=\min(\underline{\nu}-\alpha, \alpha p, \frac{\alpha k}{2})>0$ we have proved there exists a constant $C>0$ such that for all $l\in \Z$,
\begin{multline}
 \Frac{1}{2} \Frac{d}{dt} \hl^2+ m\left(\|\nabla \ul\|_{L^2}^2+ \|\nabla \ql\|_{L^2}^2+ \frac{1}{\ee^2}(\nabla \ql- \phi_\ee*\nabla \ql|\nabla \ql)_{L^2}\right)\\
\leq C(\|\nabla v\|_{L^{\infty}}+ \|v\|_{L^{\infty}}^2)\hl^2+ C\left((1+2^l)(\|\Fl\|_{L^2}+ \|\Rl\|_{L^2})+ \|\Gl\|_{L^2}+ \|\Rlp\|_{L^2}\right) \hl,
\end{multline}
which ends the proof of lemma \ref{energieZ}. $\blacksquare$\\

\subsection{Proof of the proposition}

According to the previous estimate, we observe that each term of $\hl$ except $\nabla \ql$ is parabolically regularized. The key point is that thanks to $\frac{1}{\ee^2}(\ql|\ql-\phi_\ee* \ql)_{L^2}$ (that appears in $\hl$) we will be able to get regularization on $\nabla \ql$ for some low frequencies. This leads to the hybrid Besov spaces defined in the introduction.\\

\textbf{-The low frequency case:} More precisely, the frequency threshold is simply given by the study of the decreasing function $f$ introduced in (\ref{foncaux}). As $f(z)$ goes to $1$ when $z$ goes to zero, there exists some $\gamma>0$ such that for all $0\leq z\leq \gamma$, $f(z)\geq \frac{1}{2}$. If we denote by $0<c_0<C_0$ the radii of the annulus $\cC$ used in the dyadic decomposition, let the threshold $l_\ee$ be the greatest integer such that $\ee^2 2^{2 l_\ee}C_0 \leq \gamma$ (that is $l_\ee=[\frac{1}{2}\log_2(\frac{\gamma}{C_0\ee^2})-1]$). Then for all frequency $l\leq l_\ee$ and all $\xi \in 2^l \cC$, we have:
\begin{equation}
(\ee |\xi|)^2\leq \ee^2 C_0^2 2^{2 l_\ee}\leq \gamma \quad \mbox{and} \quad 1-e^{-\ee^2 |\xi|^2}\geq \frac{\ee^2 |\xi|^2}{2},
 \label{freqseuil}
\end{equation}
and then for all $l\leq l_\ee$,
$$
(\nabla \ql- \phi_\ee*\nabla \ql|\nabla \ql)_{L^2}=C \int_{2^l \cC} (1-e^{-\ee^2 |\xi|^2}) |\hat{\nabla \ql}(\xi)|^2 d\xi \geq C\frac{\ee^2 c_0 ^2 2^{2l}}{2} \|\nabla \ql\|_{L^2}^2.
$$
So plugging this into the estimate given by lemma \ref{energieZ} gives:
\begin{multline}
 \Frac{1}{2} \Frac{d}{dt} \hl^2+ m\left(\|\nabla \ul\|_{L^2}^2+ \|\nabla \ql\|_{L^2}^2+ \frac{2^{2l}}{2\ee^2}(\ql- \phi_\ee*\ql|\ql)_{L^2}+ C\frac{c_0 ^2 2^{2l}}{2} \|\nabla \ql\|_{L^2}^2 \right)\\
\leq C(\|\nabla v\|_{L^{\infty}}+ \|v\|_{L^{\infty}}^2)\hl^2+ C\left((1+2^l)(\|\Fl\|_{L^2}+ \|\Rl\|_{L^2})+ \|\Gl\|_{L^2}+ \|\Rlp\|_{L^2}\right) \hl,
\end{multline}
and we have obtained that for all $l\leq l_\ee$:
\begin{multline}
 \Frac{d}{dt} \hl+ \frac{m}{2}2^{2l} \hl\\
\leq C(\|\nabla v\|_{L^{\infty}}+ \|v\|_{L^{\infty}}^2)\hl+ C\left((1+2^l)(\|\Fl\|_{L^2}+ \|\Rl\|_{L^2})+ \|\Gl\|_{L^2}+ \|\Rlp\|_{L^2}\right),
\end{multline}
that is for all $l\leq l_\ee$ and all $t$:
\begin{multline}
 \hl(t)+ m2^{2l}\int_0^t h_l(\tau) d\tau\\
\leq \hl(0)+ C\int_0^t (\|\nabla v\|_{L^{\infty}}+ \|v\|_{L^{\infty}}^2)\hl+ C\left((1+2^l)(\|\Fl\|_{L^2}+ \|\Rl\|_{L^2})+ \|\Gl\|_{L^2}+ \|\Rlp\|_{L^2}\right) d\tau.
\label{estimBF}
\end{multline}

\begin{rem}\sl{
 Notice that as $l_\ee \sim \log_2 (\frac{\sqrt{\gamma}}{2 C_0 \ee})$, and goes to infinity when $\ee$ goes to zero. The parabolic regularization eventually reaches every frequency.}
\end{rem}

\textbf{-The high frequency case:} Remember that $l_\ee$ is characterized by:
$$
\begin{cases}
 \ee 2^{l_\ee}C_0 \leq \sqrt{\gamma},\\
\ee 2^{l_\ee+1}C_0 > \sqrt{\gamma}
\end{cases},
$$
So, for all $l\geq l_\ee+1$ and all $\xi \in 2^l \cC$, we have
\begin{equation}
 \ee^2 |\xi|^2 \geq \ee^2 2^{2l} c_0^2 \geq \ee^2 2^{2(l_\ee+1)}C_0^2 (\frac{c_0}{C_0})^2> \gamma (\frac{c_0}{C_0})^2.
\label{borninf}
\end{equation}

Let us go back to estimate \ref{energieZ}. We can rewrite it into:
\begin{multline}
\frac{1}{2} \frac{d}{dt} \hl^2+ m 2^{2l} \left(\|\ul\|_{L^2}^2 + \|\ql\|_{L^2}^2\right) +\frac{2^{2l}}{2\ee^2} (\ql|\ql- \phi_\ee*\ql)_{L^2}+ \frac{1}{2\ee^2} (\nabla\ql|\nabla\ql- \phi_\ee* \nabla\ql)_{L^2}\\
\leq C (\|\nabla v\|_{L^\infty}+ \|v\|_{L^\infty}^2)\hl^2 + C\left((1+2^l)(\|\Fl\|_{L^2}+ \|\Rl\|_{L^2})+ \|\Gl\|_{L^2}+ \|\Rlp\|_{L^2}\right) \hl.
\end{multline}
Using Plancherel, (\ref{borninf}) together with the fact that function $x\mapsto 1-e^{-x}$ is increasing, we get that:
$$
(\nabla\ql|\nabla\ql- \phi_\ee* \nabla\ql)_{L^2}=C \int_{2^l \cC} (1-e^{-\ee^2 |\xi|^2})|\hat{\nabla \ql}(\xi)|^2 d\xi \geq (1-e^{-\gamma (\frac{c_0}{C_0})^2}) \|\nabla \ql\|_{L^2}^2.
$$
As $2^{2l}\geq \frac{\gamma}{\ee^2 C_0^2}$, (and greater than $1$ if $\ee$ is small enough) if we denote by $m'=m.\min(\frac{\gamma}{C_0^2}, 1-e^{-\gamma (\frac{c_0}{C_0})^2})>0$, our estimate becomes: for all $l\geq l_\ee+1$,
\begin{multline}
  \Frac{d}{dt} \hl+ \frac{m'}{\ee^2} \hl\\
\leq C(\|\nabla v\|_{L^{\infty}}+ \|v\|_{L^{\infty}}^2)\hl+ C\left((1+2^l)(\|\Fl\|_{L^2}+ \|\Rl\|_{L^2})+ \|\Gl\|_{L^2}+ \|\Rlp\|_{L^2}\right),
\end{multline}
that is for all $l\geq l_\ee+1$ and all $t$:
\begin{multline}
 \hl(t)+ \frac{m'}{\ee^2}\int_0^t h_l(\tau) d\tau\\
\leq \hl(0)+ C\int_0^t (\|\nabla v\|_{L^{\infty}}+ \|v\|_{L^{\infty}}^2)\hl+ C\left((1+2^l)(\|\Fl\|_{L^2}+ \|\Rl\|_{L^2})+ \|\Gl\|_{L^2}+ \|\Rlp\|_{L^2}\right) d\tau.
\label{avantparab}
\end{multline}
We need to find a way to get parabolic regularization for $\ul$. Classically (we refer for example to \cite{Dinv}, \cite{Dbook} or \cite{Has1}) we will use the new information given by the previous estimate in the energy estimate for $\ul$:
\begin{multline}
 \Frac{1}{2} \Frac{d}{dt} \|\ul\|_{L^2}^2+ \underline{\nu} 2^{2l} \|\ul\|_{L^2}^2\\
\leq C\|\nabla v\|_{L^\infty} \|\ul\|_{L^2}^2+ (\|\Gl\|_{L^2}+\|\Rlp\|_{L^2})\|\ul\|_{L^2}+ p\|\nabla \ql\|_{L^2}.\|\ul\|_{L^2}+ \frac{k}{\ee^2} \|\nabla\ql- \phi_\ee* \nabla\ql\|_{L^2}.\|\ul\|_{L^2},
\end{multline}
and then, simplifying and integrating in time:
\begin{multline}
 \|\ul(t)\|_{L^2}+ \underline{\nu} 2^{2l} \int_0^t \|\ul(\tau)\|_{L^2} d\tau \leq \|\ul(0)\|_{L^2}\\
 + C \int_0^t \left(\|\nabla v\|_{L^\infty} \|\ul\|_{L^2}+ (\|\Gl\|_{L^2}+\|\Rlp\|_{L^2}+ (p+\frac{k}{\ee^2})\|\nabla \ql\|_{L^2}\right) d\tau.
\end{multline}
Thanks to (\ref{hl}), estimate (\ref{avantparab}) implies
$$\frac{m'}{\ee^2}\int_0^t \|\nabla \ql(\tau)\|_{L^2} d\tau \leq C' \frac{m'}{\ee^2}\int_0^t h_l(\tau) d\tau,$$
so that we can estimate the last term of the left-hand side of the previous estimate and then there exists a constant such that for all $l\geq l_\ee+1$ and all $t$:
\begin{multline}
 \hl(t)+ \underline{\nu} 2^{2l} \int_0^t \|\ul(\tau)\|_{L^2} d\tau +\frac{m'}{\ee^2}\int_0^t h_l(\tau) d\tau\\
\leq \hl(0)+ C\int_0^t (\|\nabla v\|_{L^{\infty}}+ \|v\|_{L^{\infty}}^2)\hl+ C\left((1+2^l)(\|\Fl\|_{L^2}+ \|\Rl\|_{L^2})+ \|\Gl\|_{L^2}+ \|\Rlp\|_{L^2}\right) d\tau.
\label{estimHF}
\end{multline}
Remember that:
$$h_l \sim \|\ul\|_{L^2}+ \|\ql\|_{L^2} + \|\nabla\ql\|_{L^2}.$$
Let us collect what we have obtained for the low and high frequencies: if we use the Gronwall lemma to (\ref{estimBF}) and (\ref{estimHF}) we obtain that for all $t$, for all $l\leq l_\ee$:
\begin{multline}
 \|\ul(t)\|_{L^2}+ \|\ql(t)\|_{L^2} + \|\nabla\ql(t)\|_{L^2}+ m2^{2l}\int_0^t \left(\|\ul(\tau)\|_{L^2}+ \|\ql(\tau)\|_{L^2} + \|\nabla\ql(\tau)\|_{L^2}\right) d\tau\\
\leq e^{C\int_0^t (\|\nabla v(\tau)\|_{L^{\infty}}+ \|v(\tau)\|_{L^{\infty}}^2)d\tau} \Big(\|\ul(0)\|_{L^2}+ \|\ql(0)\|_{L^2} + \|\nabla\ql(0)\|_{L^2}\\
+ C\int_0^t \big((1+2^l)(\|\Fl\|_{L^2}+ \|\Rl\|_{L^2})+ \|\Gl\|_{L^2}+ \|\Rlp\|_{L^2}\big) d\tau \Big),
\end{multline}
and for $l\geq l_\ee+1$:
\begin{multline}
 \|\ul(t)\|_{L^2}+ \|\ql(t)\|_{L^2} + \|\nabla\ql(t)\|_{L^2}+ \int_0^t \left(\underline{\nu} 2^{2l} \|\ul(\tau)\|_{L^2} + \frac{m'}{\ee^2} \big(\|\ql(\tau)\|_{L^2} + \|\nabla\ql(\tau)\|_{L^2}\big)\right) d\tau\\
\leq e^{C \int_0^t (\|\nabla v(\tau)\|_{L^{\infty}}+ \|v(\tau)\|_{L^{\infty}}^2)d\tau} \Big(\|\ul(0)\|_{L^2}+ \|\ql(0)\|_{L^2} + \|\nabla\ql(0)\|_{L^2}\\
+ C\int_0^t \big((1+2^l)(\|\Fl\|_{L^2}+ \|\Rl\|_{L^2})+ \|\Gl\|_{L^2}+ \|\Rlp\|_{L^2}\big) d\tau \Big),
\end{multline}
If we use the definition of the hybrid norm introduced in the first section of this paper, after multiplying these estimates by $2^{s-1}$ and a summation over $l\in \Z$, we obtain that there exists a constant $C>0$ such that:
\begin{multline}
 \|u\|_{\tilde{L}_t^{\infty} \dot{B}_{2,1}^{s-1}}+ \|q\|_{\tilde{L}_t^{\infty} \dot{B}_{2,1}^{s-1}}+ \|q\|_{\tilde{L}_t^{\infty} \dot{B}_{2,1}^{s}}+ \|u\|_{\tilde{L}_t^1 \dot{B}_{2,1}^{s+1}}+ \|q\|_{\tilde{L}_t^1 \dot{B}_{\ee}^{s+1,s}}+ \|q\|_{\tilde{L}_t^1 \dot{B}_{\ee}^{s+2,s}}\\
\leq C e^{C\int_0^t (\|\nabla v(\tau)\|_{L^{\infty}}+ \|v(\tau)\|_{L^{\infty}}^2)d\tau} \Big(\|\ul(0)\|_{\dot{B}_{2,1}^{s-1}}+ \|\ql(0)\|_{\dot{B}_{2,1}^{s-1}} + \|\nabla\ql(0)\|_{\dot{B}_{2,1}^{s}}\\
+ \|F\|_{\tilde{L}_t^1 \dot{B}_{2,1}^{s-1}}+ \|F\|_{\tilde{L}_t^1 \dot{B}_{2,1}^{s}}+ \|G\|_{\tilde{L}_t^1 \dot{B}_{2,1}^{s-1}} +\int_0^t \big(\Sum_{l\in \Z} (2^{l(s-1)}+ 2^{ls})\|\Rl\|_{L^2}+ 2^{l(s-1)} \|\Rlp\|_{L^2}\big) d\tau \Big).
\label{quasifinale}
\end{multline}

Recall that $\Rl= [v.\nabla, \ddl] q$ and $\Rlp= [v.\nabla, \ddl] u$. We refer for example to \cite{Dbook} (lemma 2.96 and estimate 10.10) for the following well-known results:
\begin{prop}
 There exists a constant $C>0$ such that for all $v\in \dot{B}_{2,1}^{\fd+1}$ and $g\in \dot{B}_{2,1}^s$, if $r_l= [v.\nabla, \ddl] g$ there exists a nonnegative sequence $c\in l^1(\Z)$, with $\Sum_{l\in \Z} c_l=1$, such that for all $l\in \Z$:
$$
\|r_l\|_{L^2}\leq C c_l 2^{-ls} \|\nabla v\|_{\dot{B}_{2,1}^\fd} \|g\|_{\dot{B}_{2,1}^{s}}.
$$ 
\end{prop}
Then there exists three nonnegative sequences $c$, $c'$, and $c''\in l^1(\Z)$ (with $\Sum_{l\in \Z} c_l^{(,','')}=1$) such that for all $l\in \Z$:
 $$
\begin{cases}
 \|\Rl\|_{L^2}\leq C c_l 2^{-ls} \|\nabla v\|_{\dot{B}_{2,1}^\fd} \|q\|_{\dot{B}_{2,1}^{s}},\\
 \|\Rl\|_{L^2}\leq C c_l^{'} 2^{-l(s-1)} \|\nabla v\|_{\dot{B}_{2,1}^\fd} \|q\|_{\dot{B}_{2,1}^{s-1}},\\
 \|\Rlp\|_{L^2}\leq C c_l^{''} 2^{-l(s-1)} \|\nabla v\|_{\dot{B}_{2,1}^\fd} \|u\|_{\dot{B}_{2,1}^{s-1}}.
\end{cases}
$$
This allows to write:
\begin{multline}
 \int_0^t \big(\Sum_{l\in \Z} (2^{l(s-1)}+ 2^{ls})\|\Rl\|_{L^2}+ 2^{l(s-1)} \|\Rlp\|_{L^2}\big) d\tau\\
 \leq C \int_0^t \|\nabla v(\tau)\|_{\dot{B}_{2,1}^\fd} (\|q\|_{\dot{B}_{2,1}^{s-1}}+ \|q\|_{\dot{B}_{2,1}^{s}}+
\|u\|_{\dot{B}_{2,1}^{s-1}}) d\tau.
\end{multline}
Plugging this into estimate (\ref{quasifinale}), and using once more the Gronwall lemma finally gives:
 \begin{multline}
 \|u\|_{\tilde{L}_t^{\infty} \dot{B}_{2,1}^{s-1}}+ \|q\|_{\tilde{L}_t^{\infty} \dot{B}_{2,1}^{s-1}}+ \|q\|_{\tilde{L}_t^{\infty} \dot{B}_{2,1}^{s}}+ \|u\|_{\tilde{L}_t^1 \dot{B}_{2,1}^{s+1}}+ \|q\|_{\tilde{L}_t^1 \dot{B}_{\ee}^{s+1,s}}+ \|q\|_{\tilde{L}_t^1 \dot{B}_{\ee}^{s+2,s}}\\
\leq C e^{C\int_0^t (\|\nabla v(\tau)\|_{L^{\infty}}+ \|v(\tau)\|_{L^{\infty}}^2+ \|\nabla v(\tau)\|_{\dot{B}_{2,1}^\fd})d\tau} \Big(\|u_0\|_{\dot{B}_{2,1}^{s-1}}\\
+ \|q_0\|_{\dot{B}_{2,1}^{s-1}} + \|\nabla q_0\|_{\dot{B}_{2,1}^{s-1}}+ \|F\|_{\tilde{L}_t^1 \dot{B}_{2,1}^{s-1}}+ \|F\|_{\tilde{L}_t^1 \dot{B}_{2,1}^{s}}+ \|G\|_{\tilde{L}_t^1 \dot{B}_{2,1}^{s-1}}\Big).
\end{multline}
which, thanks to the injection $\dot{B}_{2,1}^\fd \hookrightarrow L^{\infty}$, ends the proof of the proposition. $\blacksquare$

\section{Proof of theorem \ref{thexist}}

The proof of the global well-posedness of system $(NSR_\ee)$ is classical and follows the lines of the proof for the compressible Navier-Stokes system (we refer to \cite{DD}, \cite{Dinv}, \cite{Dbook} (section 10.2.3) or \cite{Has1}). The only difference here is that we use, when proving the existence, the estimate given by Proposition \ref{apriori}, which ensures uniform estimates with respect to $\ee$.

The uniqueness is also classical: when $d\geq3$ we also use estimation (\ref{estimapriori}), but when $d=2$ (as for the Navier-Stokes system) we reach an endpoint for the estimates of the remainders in the paradecomposition, and the simplest is to use the same methods as for example in \cite{Dbook} (section 10.2.4 and the use of logarithmic estimates). The negative powers of $\ee$ won't cause any trouble as the uniqueness is proved when $\ee>0$ is fixed.

\subsection{Existence}

\subsubsection{Step 1: Friedrichs approximation}

In order to use the classical Friedrichs approximation, we define the frequency truncation operator $J_n$ by: for all $n\in \N$ and for all $g\in L^2(\R^d)$,
$$
J_n g= \mathcal{F}^{-1}\left(\textbf{1}_{\frac{1}{n}\leq |\xi|\leq n}(\xi) \hat{g}(\xi)\right),
$$
and we define the following approximated system (we omit the dependency in $\ee$ for more simplicity):
$$
\begin{cases}
\begin{aligned}
&\d_t q_n+ J_n\left(J_n u_n.\nabla J_n q_n\right)+ J_n\div u_n= F_n,\\
&\d_t u_n+ J_n\left(J_n u_n.\nabla J_n u_n\right) -\cA J_n u_n+P'(\overline{\rho}).\nabla J_n q_n-\frac{\kappa\overline{\rho}}{\ee^2}\nabla(\phi_\ee*J_n q_n-J_n q_n)= G_n,\\
\end{aligned}
\end{cases}
$$
where
$$
\begin{cases}
 F_n=-J_n\left(J_n q_n.\div J_n u_n\right)\\
G_n= J_n\left(K(J_n q_n).\nabla J_n q_n\right)- J_n\left(I(J_n q_n) \cA J_n u_n\right)
\end{cases}
$$
It is easy to check that it is an ordinary differential equation in $L_n^2\times(L_n^2)^d$, where $L_n^2=\{u\in L^2(\R^d), \quad J_n u=u\}$. Then for every $n\in \N$, there exists a unique maximal solution in the space $\cC^1([0, T_n^*[, L_n^2)$ and this system can be rewriten into:
$$
\begin{cases}
\begin{aligned}
&\d_t q_n+ J_n\left(u_n.\nabla q_n\right)+ \div u_n= F_n,\\
&\d_t u_n+ J_n\left(u_n.\nabla u_n\right) -\cA u_n+P'(\overline{\rho}).\nabla q_n-\frac{\kappa\overline{\rho}}{\ee^2}\nabla(\phi_\ee*q_n-q_n)= G_n,\\
\end{aligned}
\end{cases}
\leqno{(R_{\ee}^n)}
$$
where
\begin{equation}
 \begin{cases}
 F_n=-J_n\left(q_n.\div u_n\right)\\
G_n= J_n\left(K(q_n).\nabla q_n\right)- J_n\left(I(q_n) \cA u_n\right)
\end{cases}
\label{Fapprox}
\end{equation}

\subsubsection{Step 2: Uniform estimates}

Using estimate (\ref{estimapriori}) with $s=\fd$ from Proposition \ref{apriori} (as the proof is based on energy methods, the operator $J_n$ does not have any effect) we obtain that for all $t\in [0, T_n^*[$:
\begin{multline}
h_n(t) \leq C e^{C\int_0^t (\|\nabla u_n(\tau)\|_{\dot{B}_{2,1}^\fd}+ \|u_n(\tau)\|_{\dot{B}_{2,1}^\fd}^2)d\tau} \Big(\|u_n(0)\|_{\dot{B}_{2,1}^{\fd-1}}+ \|q_n(0)\|_{\dot{B}_{2,1}^{\fd-1}} + \|q_n(0)\|_{\dot{B}_{2,1}^{\fd}}\\
+ \|F_n\|_{\tilde{L}_t^1 \dot{B}_{2,1}^{\fd-1}}+ \|F_n\|_{\tilde{L}_t^1 \dot{B}_{2,1}^{\fd}}+ \|G_n\|_{\tilde{L}_t^1 \dot{B}_{2,1}^{\fd-1}}\Big).
\end{multline}
where we denote
\begin{multline}
 h_n(t)\overset{def}{=} \|u_n\|_{\tilde{L}_t^{\infty} \dot{B}_{2,1}^{\fd-1}}+ \|q_n\|_{\tilde{L}_t^{\infty} \dot{B}_{2,1}^{\fd-1}}+ \|q_n\|_{\tilde{L}_t^{\infty} \dot{B}_{2,1}^{\fd}}\\
+ \|u_n\|_{\tilde{L}_t^1 \dot{B}_{2,1}^{\fd+1}}+ \|q_n\|_{\tilde{L}_t^1 \dot{B}_{\ee}^{\fd+1,\fd}}+ \|q_n\|_{\tilde{L}_t^1 \dot{B}_{\ee}^{\fd+2,\fd}}
\end{multline}
Let $\eta>0$ be small (it will be fixed later) and assume that:
$$
h(0)\overset{def}{=} \|u(0)\|_{\dot{B}_{2,1}^{\fd-1}}+ \|q(0)\|_{\dot{B}_{2,1}^{\fd-1}} + \|q(0)\|_{\dot{B}_{2,1}^{\fd}} \leq \eta.
$$
Let us now define
$$
T_n=\sup \{t\in [0, T_n^*[,\quad h_n(t)\leq 2 C h(0)\}.
$$
As for all $n\in \N$, $h_n(0) \leq h(0) \leq \eta$, we have $T_n>0$ ($C>1$) and we will now prove by contradiction that $T_n=T_n^*$. Assume that $T_n<T_n^*$, then for all $t\leq T_n$, $h_n(t) \leq 2 C h(0) \leq 2 C \eta$ and:
$$
\int_0^t \|\nabla u_n\|_{\dot{B}_{2,1}^{\fd}} d\tau \leq h_n(t)\leq 2 C \eta
$$
and by interpolation
$$
\int_0^t \|u_n\|_{\dot{B}_{2,1}^{\fd}}^2 d\tau\leq \|u_n\|_{L^{\infty} \dot{B}_{2,1}^{\fd-1}} \|u_n\|_{L^1 \dot{B}_{2,1}^{\fd+1}} \leq h_n(t)^2\leq 4C^2\eta^2.
$$
We obtain that for all $t\leq T_n$,
$$
h_n(t) \leq Ce^{2C^2 \eta(1+2 C \eta)}\left(h_n(0)+ \|F_n\|_{L^1 \dot{B}_{2,1}^{\fd-1}} +\|F_n\|_{L^1 \dot{B}_{2,1}^{\fd}} +\|G_n\|_{L^1 \dot{B}_{2,1}^{\fd-1}}\right).
$$
We need now to estimate the right-hand side (see \ref{Fapprox} for the expressions). Thanks to estimates (\ref{estimbesov}) and Proposition \ref{estimcompo} we obtain that:
$$
\|F_n\|_{\dot{B}_{2,1}^{\fd-1}} \leq C\|u_n\|_{\dot{B}_{2,1}^{\fd+1}} \|q_n\|_{\dot{B}_{2,1}^{\fd-1}} \quad \mbox{and} \quad \|F_n\|_{\dot{B}_{2,1}^\fd} \leq C\|u_n\|_{\dot{B}_{2,1}^{\fd+1}} \|q_n\|_{\dot{B}_{2,1}^\fd}
$$
$$
\|K(q_n)\nabla q_n\|_{\dot{B}_{2,1}^{\fd-1}}\leq C\|K(q_n)\|_{\dot{B}_{2,1}^\fd} \|\nabla q_n\|_{\dot{B}_{2,1}^{\fd-1}}\leq C_\eta \|q_n\|_{\dot{B}_{2,1}^\fd}^2,
$$
$$
\|I(q_n)\cA u_n\|_{\dot{B}_{2,1}^{\fd-1}}\leq C\|I(q_n)\|_{\dot{B}_{2,1}^\fd} \|u_n\|_{\dot{B}_{2,1}^{\fd+1}}\leq C_\eta \|q_n\|_{\dot{B}_{2,1}^\fd} \|u_n\|_{\dot{B}_{2,1}^{\fd+1}}.
$$
So we can write that for all $t\leq T_n$,
\begin{multline}
h_n(t) \leq Ce^{2C^2 \eta(1+2 C \eta)}\left(h_n(0)+ \int_0^t \|u_n\|_{\dot{B}_{2,1}^{\fd+1}} (\|q_n\|_{\dot{B}_{2,1}^{\fd-1}}+ \|q_n\|_{\dot{B}_{2,1}^\fd}) d\tau\right)\\
+\int_0^t \|q_n\|_{\dot{B}_{2,1}^\fd}^2 d\tau.
\end{multline}
Using proposition \ref{estimhyb1} (see appendix) we obtain that
$$
\|q_n\|_{\dot{B}_{2,1}^\fd}^2 \leq C (\|q_n\|_{\dot{B}_{2,1}^{\fd-1}} +\|q_n\|_{\dot{B}_{2,1}^{\fd}}) \|q_n\|_{\dot{B}_{\ee}^{\fd+1,\fd}}
$$
which gives
$$
h_n(t) \leq Ce^{2C^2 \eta(1+2 C \eta)}\left(h(0)+ \int_0^t (\|u_n\|_{\dot{B}_{2,1}^{\fd+1}}+ \|q_n\|_{\dot{B}_{\ee}^{\fd+1,\fd}}) h_n(\tau) d\tau\right).
$$
Then, thanks to the Gronwall lemma and to the fact that
$$
\int_0^t (\|u_n\|_{\dot{B}_{2,1}^{\fd+1}}+ \|q_n\|_{\dot{B}_{\ee}^{\fd+1,\fd}})d\tau \leq h_n(t)\leq 2 C h(0) \leq 2 C \eta,
$$
we obtain that
$$
h_n(t) \leq Ce^{2C^2 \eta(1+2 C \eta)} h(0) e^{\int_0^t (\|u_n\|_{\dot{B}_{2,1}^{\fd+1}}+ \|q_n\|_{\dot{B}_{\ee}^{\fd+1,\fd}})d\tau} \leq C f(\eta) h(0),
$$
where
$$f(\eta)\overset{def}{=} e^{2C^2 \eta(1+2 C \eta)} e^{4C^2\eta e^{2C^2 \eta(1+2 C \eta)}}.$$
This function goes to $1$ as $\eta$ goes to zero so if $\eta>0$ is fixed and chosen so that $f(\eta)\leq \frac{3}{2}$ then for all $t\leq T_n$, $h_n(t)\leq 3Ch(0)/2< 2C h(0)$ which contradicts the fact that $T_n$ is maximal. Then $T_n=T_n^*$ and for all $t\leq T_n^*$, we have $h_n(t)\leq 2Ch(0)\leq 2C\eta$ and in particular
$$
\begin{cases}
 \|u_n\|_{L_{T_n}^\infty \dot{B}_{2,1}^{\fd-1}}<\infty\\
\|q_n\|_{L_{T_n}^\infty \dot{B}_{2,1}^{\fd-1}}<\infty
\end{cases}
$$
And, as $J_n$ is a frequency cut-off operator on $\cC(0,\frac{1}{n},n)$ we have the following bounds:
$$
\begin{cases}
 \|u_n\|_{L_{T_n}^\infty L^2}<\infty\\
\|q_n\|_{L_{T_n}^\infty L^2}<\infty
\end{cases}
$$
These bounds blow-up when $n$ goes to infinity but all that is important is that for a fixed $n$, it implies (using Cauchy-Lipschitz and the fact that $T_n^*$ is maximal) that $T_n^*= \infty$. So for all $n$ the solution $(q_n, u_n)$ is global and we have the uniform bound (in $n$ and $\ee$): 
\begin{equation}
 \forall t\geq 0, \quad h_n(t)\leq 2 Ch(0)\leq 2C \eta
\label{estimunif}
\end{equation}
Once it is obtained the rest of the method is very classical. We will use the uniform estimates to bound the time derivatives of the approximated solutions. Here will appear negative powers of $\ee$, but as it is fixed, and as only $n$ goes to infinity, the computations are strictly the same as in \cite{Dbook} or \cite{Dinv}. Then from this we easily prove that $(q_n, u_n)$ is uniformly bounded in $\cC^{\frac{1}{2}(\R_+, \dot{B}_{2,1}^{\fd-1})}\times \cC^{\frac{1}{4}(\R_+, \dot{B}_{2,1}^{\fd-\frac{3}{2}})}$ (of course, uniformly in $n$ but not in $\ee$, but $\ee$ is fixed for now). The Ascoli theorem together with the classical diagonal extraction and weak convergence arguments give, for each $\ee>0$, the existence of a global solution $(\qe, \ue)$ whose norm is bounded by $2C\eta$:
\begin{thm}
\sl{There exists $\eta>0$ such that under the same assumptions, if
$$\|u(0)\|_{\dot{B}_{2,1}^{\fd-1}}+ \|\rho(0)-\overline{\rho}\|_{\dot{B}_{2,1}^{\fd-1}} + \|\rho(0)-\overline{\rho}\|_{\dot{B}_{2,1}^{\fd}} \leq \eta$$
System $(NSR_\ee)$ has a global solution $(\rho,u)$ such that if $q=\rho-\overline{\rho}$ then $(q, u)$ is a solution in $E_\ee^s$ of $(R_{\ee})$ and
$$
\|(q,u)\|_{E_{\ee}^\fd} \leq 2 C (\|q(0)\|_{\dot{B}_{2,1}^{\fd-1} \cap \dot{B}_{2,1}^{\fd}} +\|u(0)\|_{\dot{B}_{2,1}^{\fd-1}}) \leq 2C\eta.
$$}
\end{thm}

\subsection{Uniqueness}

Once again, as $\ee$ is fixed, the computations are the same as for the compressible Navier-Stokes system. As in \cite{Dbook} or \cite{Has1}, we successively use the transport estimate for the density fluctuation (\cite{Dbook} theorem $3.14$) and the transport-diffusion estimate for the velocity (\cite{Dbook} Proposition $10.3$). If $d\geq 3$, it is simpler to use the a priori estimates provided in the present paper. We need to separate the cases $d\geq 3$ and $d=2$ (the case $d=2$ is more difficult because of endpoints for the remainder estimates in the Littlewood-Paley theory) and we obtain the following result:
\begin{thm}
 \sl{Let $d\geq 2$ and $(q_i, u_i)$ ($i\in\{1,2\}$) two solutions of $(R_\ee)$ with the same initial data on the same interval $[0, T^*[$ and both belonging to $E_\ee^\fd(T^*)$. Then there exists $\alpha>0$ such that if, for $i\in\{1,2\}$
$$
\|q_i\|_{\Tilde{L}_{T^*}^\infty (\dot{B}^{\fd-1}_{2,1}\cap\dot{B}^{\fd}_{2,1})}\leq \alpha
$$
then $(q_1,u_1)=(q_2,u_2)$ on $[0,T^*[$.
}
\end{thm}

\section{Proof of Theorem \ref{thcv}}

In this section we will show that the solution of $(R_\ee)$ goes to the solution of $(K)$ and give estimates of the speed of convergence as $\ee$ goes to zero. For that we will once more use Proposition \ref{apriori}.

As said in the introduction, if the initial data satisfy
$$
\|q_0\|_{\dot{B}_{2,1}^{\fd-1} \cap \dot{B}_{2,1}^{\fd}} +\|u_0\|_{\dot{B}_{2,1}^{\fd-1}}\leq \min(\eta_K, \eta_R),
$$
then systems $(K)$ and $(R_\ee)$ both have global solutions $(\qe,\ue)$ and $(q,u)$, and with the same notations as before, there exists a constant $C=C(\eta, \kappa)>0$ such that for all $t\in \R$ we have,
\begin{multline}
\|\ue\|_{\tilde{L}_t^{\infty} \dot{B}_{2,1}^{\fd-1}}+ \|\qe\|_{\tilde{L}_t^{\infty} \dot{B}_{2,1}^{\fd-1}}+ \|\qe\|_{\tilde{L}_t^{\infty} \dot{B}_{2,1}^{\fd}}
+\|\ue\|_{L_t^1 \dot{B}_{2,1}^{\fd+1}}\\
+\|\qe\|_{L_t^1 \dot{B}_{\ee}^{\fd+1,\fd}}+ \|\qe\|_{L_t^1 \dot{B}_{\ee}^{\fd+2,\fd}}\leq 2 C\eta,
\label{estimr}
\end{multline}
and
\begin{multline}
\|u\|_{\tilde{L}_t^{\infty} \dot{B}_{2,1}^{\fd-1}}+ \|q\|_{\tilde{L}_t^{\infty} \dot{B}_{2,1}^{\fd-1}}+ \|q\|_{\tilde{L}_t^{\infty} \dot{B}_{2,1}^{\fd}}
+\|u\|_{L_t^1 \dot{B}_{2,1}^{\fd+1}}\\
+\|q\|_{L_t^1 \dot{B}_{2,1}^{\fd+1}}+ \|q\|_{L_t^1 \dot{B}_{2,1}^{\fd+2}}\leq 2 C\eta.
\label{estimk}
\end{multline}
Up to an additional forcing term, let us write system $(K)$ with a capillary term as in system $(R_\ee)$:
$$
\begin{cases}
\begin{aligned}
&\d_t \qe+ \ue.\nabla \qe+ (1+\qe)\div \ue=0,\\
&\d_t \ue+ \ue.\nabla \ue -\cA \ue+P'(1).\nabla \qe-\frac{\kappa}{\ee^2}\nabla(\phi_\ee*\qe-\qe)= K(\qe).\nabla \qe- I(\qe) \cA \ue,\\
\end{aligned}
\end{cases}
\leqno{(R_\ee)}
$$
and
$$
\begin{cases}
\begin{aligned}
&\d_t q+ u.\nabla q+ (1+q)\div u=0,\\
&\d_t u+ u.\nabla u -\cA u+P'(1).\nabla q-\frac{\kappa}{\ee^2}\nabla(\phi_\ee*q-q)= K(q).\nabla q- I(q) \cA u+R_\ee,\\
\end{aligned}
\end{cases}
\leqno{(K)}
$$
where $R_\ee\overset{def}{=}-\frac{\kappa}{\ee^2}\nabla(\phi_\ee*q-q-\ee^2 \Delta q)$ and we recall that $K$ and $I$ are the following real-valued functions defined on $\R_+$ (whose value is $0$ when $q=0$):
$$
K(q)=\left(P'(1)-\frac{P'(1+q)}{1+q}\right) \quad \mbox{and} \quad I(q)=\frac{q}{q+1}.
$$
In order to estimate the new term $R_\ee$ we use the following lemma (which can be proved by a simple function study):
\begin{lem}
 \sl{For all $1<\beta<2$ there exists a constant $C_\beta>0$ such that for all $x\geq 0$ we have:
$$
0\leq \frac{e^{-x}-1+x}{x^\beta}\leq C_\beta
$$}
\end{lem}
so that we immediately obtain
\begin{cor}
 \sl{For all $s\in\R$ and $\beta>0$, if $q\in\dot{B}_{2,1}^{s+1+2\beta}$ then $R_\ee\in \dot{B}_{2,1}^s$ and there exists a constant $C_\beta>0$ such that we have
$$
\|R_\ee\|_{\dot{B}_{2,1}^s}\leq \kappa C_\beta \ee^{2(\beta-1)} \|q\|_{\dot{B}_{2,1}^{s+1+2\beta}}.
$$}
\label{Re}
\end{cor}
In the following we will use this estimate for $\beta$ and $s$ such that $\beta>1$ and $s+2\beta \in [\fd,\fd+1]$.
\\
Let us now write the system satisfied by the difference $(\dq,\du)=(\qe-q, \ue-u)$:
\begin{equation}
 \begin{cases}
\begin{aligned}
&\d_t \dq+ \ue.\nabla \dq+ \div \du=\delta F,\\
&\d_t \du+ \ue.\nabla \du -\cA \du+P'(1).\nabla \dq-\frac{\kappa}{\ee^2}\nabla(\phi_\ee*\dq-\dq)=\delta G-R_\ee,\\
\end{aligned}
\end{cases}
\end{equation}
where
$$
\begin{cases}
\delta F\overset{def}{=}\sum_{i=1}^3 \delta F_i\\
\delta G\overset{def}{=}\sum_{i=1}^5 \delta G_i
\end{cases}
\mbox{with }
\begin{cases}
 \delta F_1= -\du. \nabla q\\
\delta F_2= -\dq .\div \ue\\
\delta F_3= -q.\div \du
\end{cases}
\mbox{and }
\begin{cases}
 \delta G_1= -\du. \nabla u\\
 \delta G_2= \left(K(\qe)-K(q)\right).\nabla \qe\\
 \delta G_3= K(q).\nabla \dq\\
 \delta G_4= \left(I(\qe)-I(q)\right)\cA \ue\\
 \delta G_5= -I(q)\cA \du.\\
\end{cases}
$$
We will estimate $(\dq, \du)$ in less regular spaces. Let us show that for $\alpha<\min(1,d-1)$ ($\alpha$ will be precised later):
$$
\begin{cases}
 \dq \in \tilde{L}_t^{\infty} \dot{B}_{2,1}^{\fd-\alpha-1}\cap \tilde{L}_t^{\infty} \dot{B}_{2,1}^{\fd-\alpha}\cap L_t^1 \dot{B}_{\ee}^{\fd-\alpha+1,\fd-\alpha}\cap L_t^1 \dot{B}_{\ee}^{\fd-\alpha+2,\fd-\alpha}\\
\du \in \tilde{L}_t^{\infty} \dot{B}_{2,1}^{\fd-\alpha-1}\cap L_t^1 \dot{B}_{2,1}^{\fd-\alpha+1}.
\end{cases}
$$
Like previously, from estimates (\ref{estimk}) and (\ref{estimr}) we obtain (using Proposition \ref{estimhyb1} with $s=\fd$ for $\qe$, and interpolation for $q$) that $\qe$ and $q$ are uniformly bounded (with respect to $t$ and $\ee$) in $\tilde{L}_t^2 \dot{B}_{2,1}^\fd$ and that for all $p\in[1,\infty]$, $\ue$ and $u$ are uniformly bounded in $\tilde{L}_t^p \dot{B}_{2,1}^{\fd+\frac{2}{p}-1}$.
\begin{enumerate}
 \item Classically we use $\d_t \qe= -\ue.\nabla \qe -\div \ue -\qe \div \ue$, and estimate the right-hand side:
\begin{itemize}
 \item Taking $p=\frac{2}{1-\alpha}$ we get $\div \ue\in \tilde{L}_t^{\frac{2}{1-\alpha}} \dot{B}_{2,1}^{\fd-\alpha-1}$.
\item As $\ue \in\tilde{L}_t^{\frac{2}{1-\alpha}} \dot{B}_{2,1}^{\fd-\alpha}$ from (\ref{estimbesov}) we get that
$$
\begin{cases}
 \|\qe \div \ue\|_{\dot{B}_{2,1}^{\fd-\alpha-1}}\leq C \|\qe\|_{\dot{B}_{2,1}^{\fd}} \|\div \ue\|_{\dot{B}_{2,1}^{\fd-\alpha-1}}\\
 \|\ue \nabla \qe\|_{\dot{B}_{2,1}^{\fd-\alpha-1}}\leq C \|\qe\|_{\dot{B}_{2,1}^{\fd}} \|\ue\|_{\dot{B}_{2,1}^{\fd-\alpha}}.
\end{cases}
$$
\end{itemize}
So that $\d_t \qe$ is bounded (uniformly with respect to $\ee$ and $t$) in $L_t^{\frac{2}{1-\alpha}} \dot{B}_{2,1}^{\fd-\alpha-1}$ and then
$$
\qe \in \cC_t^{\frac{1+\alpha}{2}} \dot{B}_{2,1}^{\fd-\alpha-1} \quad \mbox{and }\qe-\qe(0)=\qe-q_0\in L_t^\infty \dot{B}_{2,1}^{\fd-\alpha-1}
$$
uniformly in $\ee$ but not in $t$. As the same is true for $q$, we get for all $t\geq 0$
$$
\dq\in \cC_t \dot{B}_{2,1}^{\fd-\alpha-1}\cap L_t^\infty \dot{B}_{2,1}^{\fd-\alpha-1}.
$$
By interpolation, as $\dq \in L_t^\infty \dot{B}_{2,1}^\fd$, we also have $\dq\in \cC_t \dot{B}_{2,1}^{\fd-\alpha}\cap L_t^\infty \dot{B}_{2,1}^{\fd}$ for all $t\geq 0$ and $\dq \in L_t^1 \dot{B}_{\ee}^{\fd-\alpha+1,\fd-\alpha}\cap L_t^1 \dot{B}_{\ee}^{\fd-\alpha+2,\fd-\alpha}$ (at this stage we cannot get uniform bounds in $t$ nor $\ee$, and this is not a problem for fixed $t$ and $\ee$).

\begin{rem}
 \sl{As we had to estimate $\|\ue\|_{\dot{B}_{2,1}^{\fd-\alpha}}$, this requires $\fd-\alpha\in [\fd-1, \fd+1]$ that is $\alpha \in [-1,1]$.
}
\label{Condalpha1}
\end{rem}

\item Similarly we use
$$
\d_t \ue=-\ue.\nabla \ue +\cA \ue-P'(\overline{\rho}).\nabla \qe+\frac{\kappa\overline{\rho}}{\ee^2}\nabla(\phi_\ee*\qe-\qe)+ K(\qe).\nabla \qe- I(\qe) \cA \ue,
$$
\begin{itemize}
 \item Taking $p=\frac{2}{2-\alpha}$ we get $\cA \ue\in \tilde{L}_t^{\frac{2}{2-\alpha}} \dot{B}_{2,1}^{\fd-\alpha-1}$.
\end{itemize}
\begin{rem}
 \sl{As we estimate $\|\cA \ue\|_{\dot{B}_{2,1}^{\fd-\alpha-1}}$, this requires $\fd-\alpha+1\in [\fd-1, \fd+1]$ that is $\alpha \in [0,2]$.
}
\label{Condalpha2}
\end{rem}
\begin{itemize}
\item As $\nabla q_\ee \in \tilde{L}_t^{\infty} \dot{B}_{2,1}^{\fd-1} \cap \tilde{L}_t^{\infty} \dot{B}_{2,1}^{\fd-2}$ we have $\nabla q_\ee \in \tilde{L}_t^{\infty} \dot{B}_{2,1}^{\fd-\alpha-1}$.
\item Thanks to (\ref{estimbesov}) we prove that
$$
\|\ue.\nabla \ue\|_{\tilde{L}_t^{\frac{2}{2-\alpha}} \dot{B}_{2,1}^{\fd-\alpha-1}}\leq C \|\ue\|_{\tilde{L}_t^2 \dot{B}_{2,1}^\fd} \|\nabla \ue\|_{\tilde{L}_t^{\frac{2}{1-\alpha}}\dot{B}_{2,1}^{\fd-\alpha-1}}
$$
$$
\|I(\qe)\cA \ue\|_{\tilde{L}_t^{\frac{2}{2-\alpha}} \dot{B}_{2,1}^{\fd-\alpha-1}}\leq C_\eta \|\qe\|_{\tilde{L}_t^2 \dot{B}_{2,1}^\fd} \|\cA \ue\|_{\tilde{L}_t^{\frac{2}{2-\alpha}}\dot{B}_{2,1}^{\fd-\alpha-1}}
$$
$$
\|K(\qe)\nabla \qe\|_{\dot{B}_{2,1}^{\fd-\alpha-1}}\leq C \|K(\qe)\|_{\dot{B}_{2,1}^{\fd-\alpha}} \|\nabla \qe\|_{\dot{B}_{2,1}^{\fd-1}}
$$
We obtain
$$
\|K(\qe)\nabla \qe\|_{\tilde{L}_t^\infty \dot{B}_{2,1}^{\fd-\alpha-1}}\leq C \|K(\qe)\|_{\tilde{L}_t^\infty \dot{B}_{2,1}^{\fd-\alpha}} \|\nabla \qe\|_{\tilde{L}_t^\infty \dot{B}_{2,1}^{\fd-1}}.
$$
\end{itemize}

\begin{rem}
 \sl{As we used Proposition \ref{estimcompo} we need $\fd-\alpha>0$.
}
\label{Condalpha3}
\end{rem}

So $\d_t \ue$ is bounded in the space $\tilde{L}_t^\infty \dot{B}_{2,1}^{\fd-\alpha-1}+ \tilde{L}_t^{\frac{2}{2-\alpha}}\dot{B}_{2,1}^{\fd-\alpha-1}$ (uniformly in $t$ but not in $\ee$) and then
$$
\ue \in \cC_t^{\frac{\alpha}{2}} \dot{B}_{2,1}^{\fd-\alpha-1} \quad \mbox{and }\ue-\ue(0)=\ue-u_0\in L_t^\infty \dot{B}_{2,1}^{\fd-\alpha-1}
$$
Using the same argument used for $\dq$, we obtain that $\du \in \tilde{L}_t^{\infty} \dot{B}_{2,1}^{\fd-\alpha-1}\cap L_t^1 \dot{B}_{2,1}^{\fd-\alpha+1}$.
\end{enumerate}
We can now use Proposition \ref{apriori} with $s=\fd-\alpha$ (if $\alpha<\min(1,d-1)$, the following computations are true for the cases $d=2$ and $d\geq 3$), we obtain:
\begin{multline}
 \|\du\|_{\tilde{L}_t^{\infty} \dot{B}_{2,1}^{\fd-\alpha-1}}+ \|\dq\|_{\tilde{L}_t^{\infty} \dot{B}_{2,1}^{\fd-\alpha-1}}+ \|\dq\|_{\tilde{L}_t^{\infty} \dot{B}_{2,1}^{\fd-\alpha}}+ \|\du\|_{\tilde{L}_t^1 \dot{B}_{2,1}^{\fd-\alpha+1}}\\
+\|\dq\|_{\tilde{L}_t^1 \dot{B}_{\ee}^{\fd-\alpha+1,\fd-\alpha}}+ \|\dq\|_{\tilde{L}_t^1 \dot{B}_{\ee}^{\fd-\alpha+2,\fd-\alpha}}
\leq C e^{C\int_0^t (\|\nabla \ue(\tau)\|_{\dot{B}_{2,1}^\fd}+ \|\ue(\tau)\|_{\dot{B}_{2,1}^\fd}^2)d\tau}\\
\times\Big(\|\delta F\|_{\tilde{L}_t^1 \dot{B}_{2,1}^{\fd-\alpha-1}}+\|\delta F\|_{\tilde{L}_t^1 \dot{B}_{2,1}^{\fd-\alpha}}+ \|\delta G\|_{\tilde{L}_t^1 \dot{B}_{2,1}^{\fd-\alpha-1}}+\|R_\ee\|_{\tilde{L}_t^1 \dot{B}_{2,1}^{\fd-\alpha-1}}\Big).
\label{estimdiff}
\end{multline}
Thanks to (\ref{estimbesov}) we estimate the external force terms:
\begin{itemize}
 \item $\|\delta F_1\|_{\dot{B}_{2,1}^{\fd-\alpha-1}}= \|\du.\nabla q\|_{\dot{B}_{2,1}^{\fd-\alpha-1}}\leq C\|q\|_{\dot{B}_{2,1}^{\fd+1}},\|\du\|_{\dot{B}_{2,1}^{\fd-\alpha-1}}$.
 \item $\|\delta F_1\|_{\dot{B}_{2,1}^{\fd-\alpha}} \leq C\left(\|q\|_{\dot{B}_{2,1}^{\fd+1}}.\|\du\|_{\dot{B}_{2,1}^{\fd-\alpha}} +\|q\|_{\dot{B}_{2,1}^\fd}.\|\du\|_{\dot{B}_{2,1}^{\fd-\alpha+1}}\right)$,
\item $\|\delta F_2\|_{\dot{B}_{2,1}^{\fd-\alpha-1}}=\|\dq.\div \ue\|_{\dot{B}_{2,1}^{\fd-\alpha-1}}\leq C \|\dq\|_{\dot{B}_{2,1}^{\fd-\alpha-1}} \|\ue\|_{\dot{B}_{2,1}^{\fd+1}}$,
\item $\|\delta F_2\|_{\dot{B}_{2,1}^{\fd-\alpha}}\leq C \|\dq\|_{\dot{B}_{2,1}^{\fd-\alpha}} \|\ue\|_{\dot{B}_{2,1}^{\fd+1}}$,
\item $\|\delta F_3\|_{\dot{B}_{2,1}^{\fd-\alpha-1}}=\|q\div \du\|_{\dot{B}_{2,1}^{\fd-\alpha-1}}\leq C \|q\|_{\dot{B}_{2,1}^{\fd-1}} \|\du\|_{\dot{B}_{2,1}^{\fd-\alpha+1}}$,
\item $\|\delta F_3\|_{\dot{B}_{2,1}^{\fd-\alpha}} \leq C \|q\|_{\dot{B}_{2,1}^{\fd}} \|\du\|_{\dot{B}_{2,1}^{\fd-\alpha+1}}$.
\end{itemize}
Recall that from Corollary \ref{Re},
$$
\|R_\ee\|_{\dot{B}_{2,1}^{\fd-\alpha-1}}\leq \kappa C_\beta \ee^{2(\beta-1)} \|q\|_{\dot{B}_{2,1}^{\fd+2\beta-\alpha}}.
$$
And concerning the other forcing terms in the equation on $\du$:
\begin{itemize}
 \item $\|\delta G_1\|_{\dot{B}_{2,1}^{\fd-\alpha-1}}=\|\du.\nabla u\|_{\dot{B}_{2,1}^{\fd-\alpha-1}}\leq C \|\du\|_{\dot{B}_{2,1}^{\fd-\alpha-1}} \|u\|_{\dot{B}_{2,1}^{\fd+1}}$,
\item $\|\delta G_4\|_{\dot{B}_{2,1}^{\fd-\alpha-1}}=\|\left(I(\qe)-I(q)\right)\cA \ue\|_{\dot{B}_{2,1}^{\fd-\alpha-1}}\leq C \|I(\qe)-I(q)\|_{\dot{B}_{2,1}^{\fd-\alpha}} \|\ue\|_{\dot{B}_{2,1}^{\fd+1}}$,

and using Proposition \ref{estimcompo}, we can estimate the first term:
$$
\|I(\qe)-I(q)\|_{\dot{B}_{2,1}^{\fd-\alpha}}\leq C(\|\qe\|_{L^\infty}, \|q\|_{L^\infty}) \left(|I'(0)|+\|\qe\|_{\dot{B}_{2,1}^\fd}+\|q\|_{\dot{B}_{2,1}^\fd}\right)\|\qe-q\|_{\dot{B}_{2,1}^{\fd-\alpha}},
$$
so there exists $C_\eta$ (bounded with respect to $\eta$) so that
$$
\|\delta G_4\|_{\dot{B}_{2,1}^{\fd-\alpha-1}}\leq C_\eta \|\dq\|_{\dot{B}_{2,1}^{\fd-\alpha}} \|\ue\|_{\dot{B}_{2,1}^{\fd+1}}.
$$
\item $\|\delta G_5\|_{\dot{B}_{2,1}^{\fd-\alpha-1}}=\|I(q)\cA \du\|_{\dot{B}_{2,1}^{\fd-\alpha-1}}\leq C\|I(q)\|_{\dot{B}_{2,1}^\fd} \|\du\|_{\dot{B}_{2,1}^{\fd-\alpha+1}}$,

and similarly, as $I(0)=0$, $\|I(q)\|_{\dot{B}_{2,1}^\fd}\leq C_0(\|q\|_{L^\infty}) \|q\|_{\dot{B}_{2,1}^\fd}$ so that
$$
\|\delta G_5\|_{\dot{B}_{2,1}^{\fd-\alpha-1}}\leq C_\eta \|q\|_{\dot{B}_{2,1}^\fd} \|\du\|_{\dot{B}_{2,1}^{\fd-\alpha+1}}.
$$
\end{itemize}
The last two terms require more attention: as we want to get global in time estimates, we cannot afford to loose any power of $t$ so we need to use any available $L_t^1$-type estimates of $q$ and $\qe$.
\begin{itemize}
 \item If $\alpha<d-1$ then we can use Proposition \ref{estimhyb2} with $s=\fd-\alpha$ and $t=\fd-1$ (and here $s+t>0$):
\begin{multline}
 \|\delta G_2\|_{\dot{B}_{2,1}^{\fd-\alpha-1}}=\|\left(K(\qe)-K(q)\right).\nabla \qe\|_{\dot{B}_{2,1}^{\fd-\alpha-1}}\\
\leq C \left(\|K(\qe)-K(q)\|_{\dot{B}_{2,1}^{\fd-\alpha-1}} +\|K(\qe)-K(q)\|_{\dot{B}_{2,1}^{\fd-\alpha}}\right) \|\nabla \qe\|_{\dot{B}_{\ee}^{\fd,\fd-1}}.
\end{multline}
Moreover as $\alpha\in]0,d-1[$, we have $\fd-\alpha, \fd -\alpha-1 \in ]-\fd, \fd]$, and we can use Proposition \ref{estimcompo} with $K'(0)\neq 0$, 
$$
\|K(\qe)-K(q)\|_{\dot{B}_{2,1}^{\fd-\alpha-1}} +\|K(\qe)-K(q)\|_{\dot{B}_{2,1}^{\fd-\alpha}}
\leq C_\eta \left(\|\dq\|_{\dot{B}_{2,1}^{\fd-\alpha-1}} +\|\dq\|_{\dot{B}_{2,1}^{\fd-\alpha}}\right)
$$
and we obtain:
$$
\|\delta G_2\|_{\dot{B}_{2,1}^{\fd-\alpha-1}}\leq C_\eta \left(\|\dq\|_{\dot{B}_{2,1}^{\fd-\alpha-1}} +\|\dq\|_{\dot{B}_{2,1}^{\fd-\alpha}}\right) \|\qe\|_{\dot{B}_{\ee}^{\fd+1,\fd}}.
$$
\item Similarly, Proposition \ref{estimhyb2} with $s=\fd$ and $t=\fd-\alpha-1$ ($s+t>0$) implies:
\begin{multline}
 \|\delta G_3\|_{\dot{B}_{2,1}^{\fd-\alpha-1}}=\|K(q).\nabla \dq\|_{\dot{B}_{2,1}^{\fd-\alpha-1}}\\
\leq C \left(\|K(q)\|_{\dot{B}_{2,1}^{\fd-1}} +\|K(q)\|_{\dot{B}_{2,1}^\fd}\right) \|\nabla \dq\|_{\dot{B}_{\ee}^{\fd-\alpha,\fd-\alpha-1}}.
\end{multline}
One more time, as $K(0)=0$, Proposition \ref{estimcompo} allows to write:
$$
 \|\delta G_3\|_{\dot{B}_{2,1}^{\fd-\alpha-1}}\leq C_\eta \left(\|q\|_{\dot{B}_{2,1}^{\fd-1}} +\|q\|_{\dot{B}_{2,1}^\fd}\right) \|\dq\|_{\dot{B}_{\ee}^{\fd-\alpha+1,\fd-\alpha}}.
$$
\end{itemize}
\begin{rem}
 \sl{As we used Proposition \ref{estimcompo} with $s=\fd-1$, the previous estimate is valid only when $d\geq 3$.
}
\end{rem}
In the case $d=2$, as $q\in \tilde{L}_t^\infty \dot{B}_{2,1}^{\fd-1}\cap \tilde{L}_t^1 \dot{B}_{2,1}^{\fd+1}$, thanks to interpolation, we obtain that $q\in \tilde{L}_t^2 \dot{B}_{2,1}^{\fd}$.

Similarly, but instead of using interpolation we use Proposition \ref{estimhyb1}, as $\dq\in \tilde{L}_t^\infty (\dot{B}_{2,1}^{\fd-\alpha-1}\cap\dot{B}_{2,1}^{\fd-\alpha})\cap \tilde{L}_t^1 \dot{B}_{\ee}^{\fd-\alpha+1, \fd-\alpha}$, we obtain $\dq\in \tilde{L}_t^2 \dot{B}_{2,1}^{\fd-\alpha}$. So
$$
\|\delta G_3\|_{\dot{B}_{2,1}^{\fd-\alpha-1}}\leq C \|K(q)\|_{\dot{B}_{2,1}^{\fd}} \|\nabla\dq\|_{\dot{B}_{2,1}^{\fd-\alpha-1}}.
$$
As $\|K(q)\|_{\dot{B}_{2,1}^{\fd}}\leq C(\|q\|_{L^\infty})\|q\|_{\dot{B}_{2,1}^{\fd}}\leq C_\eta\|q\|_{\dot{B}_{2,1}^{\fd}}\leq C \|q\|_{\dot{B}_{2,1}^{\fd}}$, we obtain
\begin{multline}
 \|\delta G_3\|_{\dot{B}_{2,1}^{\fd-\alpha-1}}\leq C\|q\|_{\dot{B}_{2,1}^{\fd-1}}^{\frac{1}{2}} \|q\|_{\dot{B}_{2,1}^{\fd+1}}^{\frac{1}{2}} \left( \|\dq\|_{\dot{B}_{2,1}^{\fd-\alpha}}+ \|\dq\|_{\dot{B}_{2,1}^{\fd-\alpha-1}}\right)^{\frac{1}{2}} \|\dq\|_{\dot{B}_{\ee}^{\fd-\alpha+1, \fd-\alpha}}^{\frac{1}{2}}\\
\leq \frac{C^2}{2} \|q\|_{\dot{B}_{2,1}^{\fd-1}} \|\dq\|_{\dot{B}_{\ee}^{\fd-\alpha+1, \fd-\alpha}}+ \frac{1}{2} \|q\|_{\dot{B}_{2,1}^{\fd+1}} \left( \|\dq\|_{\dot{B}_{2,1}^{\fd-\alpha}}+ \|\dq\|_{\dot{B}_{2,1}^{\fd-\alpha-1}}\right).
\end{multline}

\begin{rem}
 \sl{As we used Proposition \ref{estimcompo} we need the sum of indices $s$ and $t$ to be positive so $\alpha<d-1$.
}
\label{Condalpha4}
\end{rem}
If we define
\begin{multline}
h(t)\overset{def}{=}\|\du\|_{\tilde{L}_t^{\infty} \dot{B}_{2,1}^{\fd-\alpha-1}}+ \|\dq\|_{\tilde{L}_t^{\infty} \dot{B}_{2,1}^{\fd-\alpha-1}}+ \|\dq\|_{\tilde{L}_t^{\infty} \dot{B}_{2,1}^{\fd-\alpha}}+ \|\du\|_{\tilde{L}_t^1 \dot{B}_{2,1}^{\fd-\alpha+1}}\\
+\|\dq\|_{\tilde{L}_t^1 \dot{B}_{\ee}^{\fd-\alpha+1,\fd-\alpha}}+ \|\dq\|_{\tilde{L}_t^1 \dot{B}_{\ee}^{\fd-\alpha+2,\fd-\alpha}}
\end{multline}
then collecting the previous estimates on the external terms in (\ref{estimdiff}), there exists $C_\eta>0$ (bounded when $\eta\in[0,1]$) such that we have for all $t\geq 0$ (and in any case $d=2$ or $d\geq 3$),
\begin{multline}
h(t)\leq C_\eta e^{2C\eta+ 4C^2 \eta^2} \Big(\int_0^t h(\tau) \big(\|q\|_{\dot{B}_{2,1}^{\fd+1}} +\|q\|_{\dot{B}_{2,1}^{\fd+2}} +\|\ue\|_{\dot{B}_{2,1}^{\fd+1}}\\
+\|u\|_{\dot{B}_{2,1}^{\fd+1}} +\|\qe\|_{\dot{B}_{\ee}^{\fd+1,\fd}}\big)d\tau +\kappa C_\beta \ee^{2(\beta-1)} \|q\|_{\Tilde{L}_t^1\dot{B}_{2,1}^{\fd+2\beta-\alpha}}\\
+\big(\|q\|_{\Tilde{L}_t^\infty \dot{B}_{2,1}^{\fd-1}} +\|q\|_{\Tilde{L}_t^\infty \dot{B}_{2,1}^\fd}\big).\big(\|\dq\|_{\Tilde{L}_t^1 \dot{B}_{\ee}^{\fd-\alpha+1,\fd-\alpha}} +\|\du\|_{\Tilde{L}_t^1 \dot{B}_{2,1}^{\fd+1-\alpha}} \big)\Big).
\end{multline}
Thanks to (\ref{estimr}) and (\ref{estimk}), this estimate turns into:
\begin{multline}
h(t)\leq C_\eta \Big(\int_0^t h(\tau) \big(\|q\|_{\dot{B}_{2,1}^{\fd+1}} +\|q\|_{\dot{B}_{2,1}^{\fd+2}} +\|\ue\|_{\dot{B}_{2,1}^{\fd+1}}+\|u\|_{\dot{B}_{2,1}^{\fd+1}} +\|\qe\|_{\dot{B}_{\ee}^{\fd+1,\fd}}\big)d\tau\\
+\kappa C_\beta \ee^{2(\beta-1)} \|q\|_{\Tilde{L}_t^1\dot{B}_{2,1}^{\fd+2\beta-\alpha}}+4\eta h(t)\Big).
\label{estimpresquefinale}
\end{multline}
The term $\|q\|_{\Tilde{L}_t^1\dot{B}_{2,1}^{\fd+2\beta-\alpha}}$ can be estimated by interpolation (thanks again to (\ref{estimk})) if $1\leq 2\beta-\alpha\leq 2$ and in this case it is less than $2C.\eta$. Let us collect all the conditions on $\alpha$ (see remarks \ref{Condalpha1} to \ref{Condalpha4}) and $\beta$:
$$
\begin{cases}
\alpha \in ]0,1],\quad \alpha<\fd, \quad \alpha< d-1,\\
\beta\in]1,2[,\quad \frac{1}{2}+\frac{\alpha}{2}\leq \beta\leq 1+\frac{\alpha}{2}.
\end{cases}
$$
This obviously implies that $d\geq 2$, and if $\alpha$ is fixed in $]0,1[$ (when $d=2$) or in $]0,1]$ (when $d\geq 3$), and if we choose $\beta=1+\frac{\alpha}{2}$ all the conditions are satisfied. Moreover if $\eta>0$ is such that $4 \eta. C_\eta\leq \frac{1}{2}$ then the last term in \ref{estimpresquefinale} can be absorbed by the left-hand side, and tanks to the Gronwall lemma, for all $t\geq 0$, we have
$$
h(t)\leq C_{\eta,\alpha} \kappa \ee^{\alpha} e^{\int_0^t \big(\|q\|_{\dot{B}_{2,1}^{\fd+1}} +\|q\|_{\dot{B}_{2,1}^{\fd+2}} +\|\ue\|_{\dot{B}_{2,1}^{\fd+1}}+\|u\|_{\dot{B}_{2,1}^{\fd+1}} +\|\qe\|_{\dot{B}_{\ee}^{\fd+1,\fd}}\big)d\tau}
$$
And finally, thanks to (\ref{estimk}) and (\ref{estimr}), we obtain that for all $t\in\R$,
\begin{multline}
 \|\du\|_{\tilde{L}_t^{\infty} \dot{B}_{2,1}^{\fd-\alpha-1}}+ \|\dq\|_{\tilde{L}_t^{\infty} \dot{B}_{2,1}^{\fd-\alpha-1}}+ \|\dq\|_{\tilde{L}_t^{\infty} \dot{B}_{2,1}^{\fd-\alpha}}+ \|\du\|_{\tilde{L}_t^1 \dot{B}_{2,1}^{\fd-\alpha+1}}\\
+\|\dq\|_{\tilde{L}_t^1 \dot{B}_{\ee}^{\fd-\alpha+1,\fd-\alpha}}+ \|\dq\|_{\tilde{L}_t^1 \dot{B}_{\ee}^{\fd-\alpha+2,\fd-\alpha}}\leq C_{\eta,\alpha} \kappa \ee^{\alpha}.
\end{multline}
which ends the proof of the theorem. $\blacksquare$

\begin{rem}
 If we had not estimated this way the terms $\delta G_2$ and $\delta G_3$, we would write that $\delta G_2+ \delta G_3=\nabla (L(\qe)-L(q))$ and then due to the regularity of $\dq$, we could only obtain that $\|\delta G_2+\delta G_2\|_{L_t^1\dot{B}_{2,1}^{\fd-\alpha-1}}\leq tC_\eta \|\dq\|_{L_t^\infty \dot{B}_{2,1}^{\fd-\alpha}}$, which prevents any global in time estimate.
\end{rem}

\section{Appendix}

In this section we state and prove some estimates involving the hybrid Besov spaces introduced in the beginning of the paper. These estimates are very close to the one proved in \cite{Dinv} (see appendix).
\begin{prop}
 \sl{
Let $s\in \R$, $\alpha>0$. For all $q\in \Tilde{B}_{\alpha}^{s,\infty}\cap \Tilde{B}_{\alpha}^{s,1}$, we have
$$
\|q\|_{\dot{B}_{2,1}^s}^2 \leq \|q\|_{\Tilde{B}_{\alpha}^{s,\infty}} \|q\|_{\Tilde{B}_{\alpha}^{s,1}}
$$
}
\end{prop}
\textbf{Proof :} This proof is classical (see \cite{Dinv} appendix). From the definitions given in the first section we have:
\begin{multline}
 \|q\|_{\dot{B}_{2,1}^s}^2= \Sum_{l\in\Z} 2^{ls}\|\ddl q\|_{L^2}\\
=\Sum_{l\in\Z} (2^{ls}\max(\alpha,2^{-l})\|\ddl q\|_{L^2})^{\frac{1}{2}} (\frac{2^{ls}}{\max(\alpha,2^{-l})}\|\ddl q\|_{L^2})^{\frac{1}{2}}\\
\leq \left(\Sum_{l\in\Z} 2^{ls}\max(\alpha,2^{-l})\|\ddl q\|_{L^2}\right)^{\frac{1}{2}} \left(\Sum_{l\in\Z} \frac{2^{ls}}{\max(\alpha,2^{-l})}\|\ddl q\|_{L^2}\right)^{\frac{1}{2}}\\
\leq \|q\|_{\Tilde{B}_{\alpha}^{s,\infty}} \|q\|_{\Tilde{B}_{\alpha}^{s,1}}.\blacksquare
\end{multline}
\begin{rem}
  \sl{More generally, if $s\in \R$, $\alpha>0$, $r,r'\in [1,\infty]$ with $\frac{1}{r}+\frac{1}{r'}$. For all $q\in \Tilde{B}_{\alpha}^{s,r}\cap \Tilde{B}_{\alpha}^{s,r'}$, we have
$$
\|q\|_{\dot{B}_{2,1}^s}^2 \leq \|q\|_{\Tilde{B}_{\alpha}^{s,r}} \|q\|_{\Tilde{B}_{\alpha}^{s,r'}}
$$
}
\end{rem}

\begin{prop}
\sl{Let $s\leq d/2$, $t\leq d/2-1$ such that $s+t>0$. There exists $C>0$ such that for all $(u,v)\in \Tilde{B}_{\alpha}^{s,\infty}\times \Tilde{B}_{\alpha}^{t,1}$,
$$
\|uv\|_{\dot{B}_{2,1}^{s+t-\fd}}\leq C \|u\|_{\Tilde{B}_{\alpha}^{s,\infty}} \|v\|_{\Tilde{B}_{\alpha}^{t,1}}.
$$
}
\end{prop}
\textbf{Proof :} Using the Bony decomposition, we have $uv=T_u v+ T_v u+ R(u,v)$ (we refer to the introduction for the definitions of paraproduct and remainder) and we will separately estimate these terms in $\dot{B}_{2,1}^{s+t-\fd}$. Before that let us state the following lemma (we refer to \cite{Dinv} Proposition $5.3$.)
\begin{lem}
 \sl{Let $\alpha>0$, $a,b\in\R$. Then we have
$$
\frac{\max(\alpha,2^{-a})}{\max(\alpha,2^{-b})}\leq
\begin{cases}
 1 & \mbox{if}\quad a\geq b,\\
2^{b-a} & \mbox{if}\quad a\leq b.
\end{cases}
$$}
\label{estimax}
\end{lem}
\textbf{1-} Let us begin with $T_u v= \sum_{q\in\Z} S_{q-1}u. \Delta_q v$. As for all $q\in \Z$, $S_{q-1}u. \Delta_q v$ has its frequencies localized in $2^q \cC'$ (where $\cC'$ is a ring) we will bound the norm by estimating $2^{q(s+t-\fd)}\|S_{q-1}u. \Delta_q v\|_{L^2}$. Using the Bernstein lemma,
\begin{multline}
 2^{q(s+t-\fd)}\|S_{q-1}u. \Delta_q v\|_{L^2}\leq 2^{q(s+t-\fd)} \left(\Sum_{q'\leq q-2} \|\Delta_{q'} u\|_{L^\infty}\right).\|\Delta_q v\|_{L^2}\\
\leq 2^{q(s+t-\fd)} \left(\Sum_{q'\leq q-2} 2^{q'\fd}\|\Delta_{q'} u\|_{L^2}\right).\|\Delta_q v\|_{L^2}
\end{multline}
As $(u,v)\in \Tilde{B}_{\alpha}^{s,\infty}\times \Tilde{B}_{\alpha}^{t,1}$, there exist two nonnegative sequences $c,c'\in l^1(\Z)$ such that $\|c\|_{l^1}\leq 1$, $\|c'\|_{l^1}\leq 1$ and for all $j\in \Z$,
\begin{equation}
 \|\Delta_j u\|_{L^2} \leq \frac{2^{-js} c_j}{\max(\alpha, 2^{-j})}\|u\|_{\Tilde{B}_{\alpha}^{s,\infty}} \quad \mbox{and} \quad \|\Delta_j v\|_{L^2} \leq 2^{-jt} c_j'\max(\alpha, 2^{-j})\|v\|_{\Tilde{B}_{\alpha}^{s,1}},
\label{uv}
\end{equation}
so we can write that
$$
2^{q(s+t-\fd)}\|S_{q-1}u. \Delta_q v\|_{L^2}\leq \Sum_{q'\leq q-2} 2^{(q-q')(s-\fd)} c_{q'} c_q' \frac{\max(\alpha,2^{-q})}{\max(\alpha,2^{-q'})} \|u\|_{\Tilde{B}_{\alpha}^{s,\infty}} \|v\|_{\Tilde{B}_{\alpha}^{s,1}}.
$$
Thanks to the previous lemma, as $q'\leq q-2$ we obtain
$$
2^{q(s+t-\fd)}\|S_{q-1}u. \Delta_q v\|_{L^2}\leq \left(\Sum_{q'\leq q-2} 2^{(q-q')(s-\fd)} c_{q'}\right)c_q' \|u\|_{\Tilde{B}_{\alpha}^{s,\infty}} \|v\|_{\Tilde{B}_{\alpha}^{s,1}}.
$$
If we denote
$$
a_q=
\begin{cases}
 2^{q(s-\fd)} & \mbox{if } q\geq 2\\
0 & \mbox{if } q<2,
\end{cases}
$$
then $\left(\Sum_{q'\leq q-2} 2^{(q-q')(s-\fd)} c_{q'}\right)=(a*c)_q$, and as $c,c'\in l^1(\Z)$, $a\in l^{\infty}(\Z)$ if and only if $s\leq \fd$ then $(a*c).c' \in l^1(\Z)$ when $s\leq \fd$ and
$$
\|T_u v\|_{\dot{B}_{2,1}^{s+t-\fd}}\leq C \|u\|_{\Tilde{B}_{\alpha}^{s,\infty}} \|v\|_{\Tilde{B}_{\alpha}^{s,1}}.
$$
\textbf{2-} Let us now turn to $T_v u$. The same arguments as previously give that for all $q\in \Z$,
$$
2^{q(s+t-\fd)}\|S_{q-1}v. \Delta_q u\|_{L^2}\leq \Sum_{q'\leq q-2} 2^{(q-q')(t-\fd)} c_{q'} c_q' \frac{\max(\alpha,2^{-q'})}{\max(\alpha,2^{-q})} \|u\|_{\Tilde{B}_{\alpha}^{s,\infty}} \|v\|_{\Tilde{B}_{\alpha}^{s,1}}.
$$
This time, as $q'\leq q-2$ thanks to the previous lemma, $\frac{\max(\alpha,2^{-q'})}{\max(\alpha,2^{-q})}\leq 2^{q-q'}$ and
$$
2^{q(s+t-\fd)}\|S_{q-1}u. \Delta_q v\|_{L^2}\leq \left(\Sum_{q'\leq q-2} 2^{(q-q')(t-\fd+1)} c_{q'}\right)c_q' \|u\|_{\Tilde{B}_{\alpha}^{s,\infty}} \|v\|_{\Tilde{B}_{\alpha}^{s,1}}.
$$
Similarly, if we denote
$$
b_q=
\begin{cases}
 2^{q(t-\fd+1)} & \mbox{if } q\geq 2\\
0 & \mbox{if } q<2,
\end{cases}
$$
Then if $t\leq \fd-1$ we get that $\|(b*c).c'\|_{l^1(\Z)}\leq \|b\|_{l^\infty(\Z)}\|c\|_{l^1(\Z)}\|c'\|_{l^1(\Z)}$ and
$$
\|T_v u\|_{\dot{B}_{2,1}^{s+t-\fd}}\leq C \|u\|_{\Tilde{B}_{\alpha}^{s,\infty}} \|v\|_{\Tilde{B}_{\alpha}^{s,1}}.
$$
\textbf{3-} Finally, let us look at the remainder $R(u,v)=\Delta_q u(\Delta_{q-1} v+ \Delta_q v +\Delta_{q+1} v)$. Here when $q\in\Z$ the frequencies are only in the ball $2^q \mathcal{B}$ and we estimate, for all $j$, $2^{j(s+t-\fd)} \|\Delta_j R(u,v)\|_{L^2}$. If we perform the same computations as before (using $\|uv\|_{L^2}\leq \|u\|_{L^infty} \|v\|_{L^2}$ and the Bernstein lemma to estimate the $L^\infty$-norm by the $L^2$-norm as the freuencies are localized), we will end with the condition $s+t>\fd$. We can get a better condition if we first estimate the $L^2$-norm by the $L^1$-norm (using the spectral localization) and then use the Hölder estimate. Thanks to the Bernstein lemma:
\begin{multline}
 2^{j(s+t-\fd)} \|\Delta_j R(u,v)\|_{L^2} \leq 2^{j(s+t-\fd)} 2^{j\fd} \|\Delta_j R(u,v)\|_{L^1}\\
\leq 2^{j(s+t)} \Sum_{q\geq j-N_0} \|\Delta_q u\|_{L^2} (\|\Delta_{q-1} v\|_{L^2}+ \|\Delta_q v\|_{L^2} +\|\Delta_{q+1} v\|_{L^2}).
\end{multline}
As $(u,v)\in \Tilde{B}_{\alpha}^{s,\infty}\times \Tilde{B}_{\alpha}^{t,1}$, we can use (\ref{uv}) and get:
\begin{multline}
 2^{j(s+t-\fd)} \|\Delta_j R(u,v)\|_{L^2} \leq 2^{j(s+t)} \Sum_{q\geq j-N_0} \frac{2^{-qs} c_q}{\max(\alpha, 2^{-q})}\|u\|_{\Tilde{B}_{\alpha}^{s,\infty}}\times 2^{-qt} \|v\|_{\Tilde{B}_{\alpha}^{s,1}}\\
\times \left(2^t c_{q-1}'\max(\alpha, 2^{-(q-1)}+ c_q'\max(\alpha, 2^{-q}) +2^{-t} c_{q+1}'\max(\alpha, 2^{-(q+1)}\right).
\end{multline}
Thanks again to lemma \ref{estimax} and (\ref{uv}), we obtain:
\begin{multline}
 2^{j(s+t-\fd)} \|\Delta_j R(u,v)\|_{L^2}\\
\leq (1+2^t+2^{-t}) \|u\|_{\Tilde{B}_{\alpha}^{s,\infty}} \|v\|_{\Tilde{B}_{\alpha}^{s,1}} \Sum_{q\geq j-N_0} 2^{(j-q)(s+t)} c_q(c_{q-1}'+c_q'+c_{q+1}')\\
\leq (1+2^t+2^{-t}) \|u\|_{\Tilde{B}_{\alpha}^{s,\infty}} \|v\|_{\Tilde{B}_{\alpha}^{s,1}} \Sum_{q\geq j-N_0} 2^{(j-q)(s+t)} c_q.
\end{multline}
Once again, if we denote
$$
d_q=
\begin{cases}
 2^{j(s+t)} & \mbox{if } j\leq N_0\\
0 & \mbox{else},
\end{cases}
$$
As we want $(d*c)\in l^1(\Z)$ we now need that $d\in l^1(\Z)$ ($d\in l^\infty(\Z)$ is not enough as we only could estimate $c_{q-1}'+c_q'+c_{q+1}'\leq 1$ in $L^\infty$ instead of $L^1$) which is true if and only if $s+t>0$. So when $s+t>0$ we finally get:
$$
\|R(u,v)\|_{\dot{B}_{2,1}^{s+t-\fd}}\leq C_{s,t} \|u\|_{\Tilde{B}_{\alpha}^{s,\infty}} \|v\|_{\Tilde{B}_{\alpha}^{s,1}}.
$$
which ends the proof of the proposition. $\blacksquare$

\begin{rem}
\sl{
For all $q\in \dot{B}_{2,1}^{s-1}\cap \dot{B}_{2,1}^{s}= \Tilde{B}_1^{s,\infty}$ we have
$$
\|q\|_{\Tilde{B}_1^{s,\infty}}\leq \|q\|_{\dot{B}_{2,1}^{s-1}} +\|q\|_{\dot{B}_{2,1}^{s}}
$$
and thanks to remark \ref{hybridcomparaison}, when $\ee>0$ is small enough, for all $q\in \dot{B}_{\ee}^{s+1,s}$, we have
$$
\|q\|_{\Tilde{B}_1^{s,1}}\leq \|q\|_{\dot{B}_{\ee}^{s+1,s}},
$$
so we can use the hybrid norms introduced in (\ref{normhybride}) and we will in fact use the following results:}
\end{rem}
\begin{prop}
 \sl{Let $s\in \R$. There exists a constant $C>0$ such that for all $\ee>0$, and all $q\in \dot{B}_{2,1}^{s-1}\cap \dot{B}_{2,1}^{s}\cap \dot{B}_{\ee}^{s+1,s}$, we have
$$
\|q\|_{\dot{B}_{2,1}^s}^2 \leq C (\|q\|_{\dot{B}_{2,1}^{s-1}} +\|q\|_{\dot{B}_{2,1}^{s}}) \|q\|_{\dot{B}_{\ee}^{s+1,s}}
$$}
\label{estimhyb1}
\end{prop}
\begin{prop}
\sl{Let $s\leq d/2$, $t\leq d/2-1$ such that $s+t>0$. There exists $C>0$ such that for all $\ee>0$, and all $(u,v)\in (\dot{B}_{2,1}^{s-1}\cap \dot{B}_{2,1}^{s})\times \dot{B}_{\ee}^{t+1,t}$,
$$
\|uv\|_{\dot{B}_{2,1}^{s+t-\fd}}\leq C (\|u\|_{\dot{B}_{2,1}^{s-1}} +\|u\|_{\dot{B}_{2,1}^{s}}) \|v\|_{\dot{B}_{\ee}^{t+1,t}}.
$$}
\label{estimhyb2}
\end{prop}

\end{document}